\documentclass{article}
\usepackage{amssymb,amsmath,amsfonts,amsthm}
\usepackage{graphics}
\usepackage[matrix,arrow,curve]{xy}

\newtheoremstyle{my}{1.5em}{0.5em}{\em}{}{\sc}{.}{0.5em}{}
\newtheoremstyle{mydef}{1.5em}{0.5em}{}{}{\sc}{.}{0.5em}{}

\theoremstyle{my}
\newtheorem{thm}{Theorem}[section]
\newtheorem{theorem}[thm]{Theorem}
\newtheorem{cor}[thm]{Corollary}

\newtheorem{lemma}[thm]{Lemma}

\newtheorem{prop}[thm]{Proposition}

\newtheorem{defn}[thm]{Definition}

\newtheorem{remark}[thm]{Remark}

\newtheorem{example}[thm]{Example}

\newcommand{\acknowledgments}{{\em Acknowledgments.} }

\newcommand{\R}{\mathbb{R}}
\newcommand{\Z}{\mathbb{Z}}
\newcommand{\Q}{\mathbb{Q}}
\newcommand{\C}{\mathbb{C}}

\newcommand{\half}{{\textstyle\frac{1}{2}}}

\newcommand{\iso}{\cong}           
\newcommand{\htp}{\simeq}          
\newcommand{\smooth}{C^\infty}
\newcommand{\CP}[1]{\C {\mathrm P}^{#1}}

\newcommand{\leftsc}{\langle}
\newcommand{\rightsc}{\rangle}

\newcommand{\suchthat}{\; : \;}

\renewcommand{\ker}{\mathrm{ker}}

\renewcommand{\o}{\omega}

\newcommand{\Diff}{\text{\it Diff}\,}

{

\begin{enumerate}
\itemsep1em \leftmargin2em
}{\end{enumerate}}

{\begin{list}{{\rm(\alph{enumi}) }}{\usecounter{enumi}

\leftmargin0cm \labelsep0cm \rightmargin0cm \parsep1em
\listparindent0em \itemsep-1em \topsep1em \parskip1em
\setlength{\labelwidth}{\fill}}} {\parskip0em \end{list}}

{

\begin{enumerate}}{\end{enumerate}}

{

\begin{enumerate} \parsep0cm
\itemsep0cm \parskip0cm}{\end{enumerate}}

{

\begin{enumerate} \parsep0cm
\itemsep0cm \parskip0cm}{\end{enumerate}}

{

\begin{enumerate}}{\end{enumerate}}

{

\begin{enumerate}}{\end{enumerate}}

{\begin{list}{(\alph{enumi}) }{\usecounter{enumi}

\leftmargin0cm \labelsep0cm
\itemsep1em \rightmargin0cm
\setlength{\labelwidth}{\fill}}} {\end{list}}

\newcommand{\BB}{\mathcal{B}}
\newcommand{\EE}{\mathcal{E}}
\newcommand{\LL}{\mathcal{L}}
\newcommand{\OO}{\mathcal{O}}
\renewcommand{\SS}{\mathcal{S}}

\newcommand{\MM}{\mathcal{M}}
\newcommand{\JJ}{\mathcal{J}}
\newcommand{\gen}[1]{\langle #1 \rangle}
\newcommand{\K}{\mathbb{K}}
\newcommand{\G}{\mathcal{G}}

\numberwithin{equation}{section}

\renewcommand{\subsection}{\hspace{-\parindent}\refstepcounter{subsection}{\bf
(\arabic{section}\alph{subsection}) }}

\begin{document}
\title{Lectures on \\ four-dimensional Dehn twists}
\author{Paul Seidel}
\date{September 1, 2003}
\maketitle

\setcounter{section}{-1}
\section{Introduction}

Let $M$ be a closed symplectic manifold, with symplectic form $\o$. A
symplectic automorphism is a diffeomorphism $\phi: M \rightarrow M$
such that $\phi^*\o = \o$. We equip the group $Aut(M) = Aut(M,\o)$ of
all such maps with the $\smooth$-topology. Like the whole
diffeomorphism group, this is an infinite-dimensional Lie group in a
very loose sense: it has a well-defined Lie algebra, which consists
of closed one-forms on $M$, but the exponential map is not locally
onto. We will be looking at the homotopy type of $Aut(M)$, and in
particular the {\em symplectic mapping class group} $\pi_0(Aut(M))$.

\begin{remark} \label{th:aut-h}
If $H^1(M;\R) \neq 0$, the $\smooth$-topology is is many respects not
the right one, and should be replaced by the {\em Hamiltonian
topology}, denoted by $Aut^h(M)$. This is defined by taking a basis
of neighbourhoods of the identity to be the symplectic automorphisms
generated by time-dependent Hamiltonians $H: [0;1] \times M
\rightarrow \R$ with $||H||_{C^k} < \epsilon$ for some $k,\epsilon$.
A smooth isotopy is continuous in the Hamiltonian topology iff it is
Hamiltonian. The relation between $\pi_0(Aut(M))$ and
$\pi_0(Aut^h(M))$ is determined by the image of the flux
homomorphism, which we do not discuss since it is thoroughly covered
elsewhere \cite{lalonde-mcduff-polterovich95}. In fact, for
simplicity we will mostly use $Aut(M)$, even when this restricts us
to manifolds with $H^1(M;\R) = 0$ (if this irks, see Remark
\ref{th:hamiltonian}).
\end{remark}

When $M$ is two-dimensional, Moser's lemma tells us that $\Diff^+(M)$
retracts onto $Aut(M)$, so $\pi_0(Aut(M))$ is the ordinary mapping
class group, which leaves matters in the hands of topologists. Next,
suppose that $M$ is a four-manifold. Diffeomorphism groups in four
dimensions are not well understood, not even the local case of
$\R^4$. Contrarily to what this seems to indicate, the corresponding
symplectic problem is far easier: this was one of 
Gromov's \cite{gromov85} original applications of the pseudo-holomorphic curve
method. In extreme simplification, the strategy is to fibre a symplectic four-manifold
by a family of such curves, and thereby to reduce the isotopy
question to a fibered version of the two-dimensional case. For
instance, it turns out that the compactly supported symplectic
automorphism group $Aut^c(\R^4)$ is weakly contractible. Here are some 
more of Gromov's results:

\begin{theorem} \label{th:gromov}
(1) $Aut(\CP{2})$ is homotopy equivalent to $PU(3)$. (2) For a
monotone symplectic structure, $Aut(S^2 \times S^2)$ is homotopy
equivalent to $(SO(3) \times SO(3)) \rtimes \Z/2$. (3) (not actually
stated in \cite{gromov85}, but follows by the same method) for a
monotone symplectic structure, $Aut(\CP{2} \# \overline{\CP{2}})$ is
homotopy equivalent to $U(2)$.
\end{theorem}

Recall that a symplectic manifold is monotone if $c_1(M) = r [\o]
\in H^2(M;\R)$ for some $r>0$. (Our formulation is slightly
anachronistic: it is true that symplectic forms on $\CP{2}$,
$S^2 \times S^2$, and $\CP{2} \# \overline{\CP{2}}$ are determined
up to isomorphism by their cohomology classes, but this is a more recent result, whose
proof depends on Seiberg-Witten invariants and Taubes' work; originally, Theorem
\ref{th:gromov} would have 
been formulated in terms of monotone K{\"a}hler forms, which obviously give rise
to unique symplectic structures.) Note that in all cases, the
result says that $Aut(M)$ is homotopy equivalent to the group of
holomorphic automorphisms. One can average the K{\"a}hler form
with respect to a maximal compact subgroup of this, and then
$Aut(M)$ becomes homotopy equivalent to the K{\"a}hler isometry
group.

After surmonting considerable difficulties, Abreu and McDuff
\cite{abreu97,abreu-mcduff98} (see also \cite{anjos02}) extended
Gromov's method to non-monotone symplectic forms. Their results show
that the symplectic automorphism group changes radically if one
varies the symplectic class. Moreover, it is not typically homotopy
equivalent to any compact Lie group, so that K{\"a}hler isometry
groups are no longer a good model. Nevertheless, they obtained an
essentially complete understanding of the topology of $Aut(M)$, in
particular:

\begin{theorem} \label{th:abreu-mcduff}
Suppose that $M$ is either $S^2 \times S^2$ or $\CP{2} \#
\overline{\CP{2}}$, with a non-monotone symplectic form. Then
$\pi_0(Aut(M))$ is trivial.
\end{theorem}

Now we bring a different source of intuition into play. Let $B$ be
a connected pointed manifold. A symplectic fibration with fibre
$M$ and base $B$ is a smooth proper fibration $\pi: E \rightarrow
B$ together with a family of symplectic forms $\{\Omega_b\}$ on
the fibres such that $[\Omega_b] \in H^2(E_b;\R)$ is locally
constant, and a preferred isomorphism between the fibre over the
base point and $M$. There is a universal fibration in the homotopy
theory sense, whose base is the classifying space $BAut(M)$. The
main advantage of the classifying space viewpoint is that it
provides a link with algebraic geometry. Namely, let $\EE$, $\BB$
be smooth quasi-projective varieties, and $\pi: \EE \rightarrow
\BB$ a proper smooth morphism of relative dimension two, with a
line bundle $\LL \rightarrow \EE$ which is relatively very ample.
This means that the sections of $\LL|\EE_b$ define an embedding of
$\EE$ into a projective bundle over $\BB$. From this embedding one
can get a family of (Fubini-Study) K{\"a}hler forms on the fibres,
so $\EE$ becomes a symplectic fibration, classified by a map
\begin{equation} \label{eq:bb}
\BB \longrightarrow BAut(\MM)
\end{equation}
where $\MM$ is the fibre over some base point, equipped with its
K{\"a}hler form. In some cases, one can construct a family which is
universal in the sense of moduli theory, and then the associated map
\eqref{eq:bb} is the best of its kind. More generally, one needs to
consider versal families together with the automorphism groups of
their fibres (this is very much the case in the situation studied by
Abreu and McDuff; it would be nice to have a sound stack-theoretic
formulation, giving the right generalization of the universal base
space at least as a homotopy type). Of course, there is no {\em a
priori} guarantee that algebraic geometry comes anywhere near
describing the whole topology of the symplectic automorphism group,
or vice versa, that symplectic topology detects all of the structure
of algebro-geometric moduli spaces.

\begin{example}
Suppose that some $\MM$ is a double cover of $\CP{1} \times
\CP{1}$ branched along a smooth curve of bidegree (6,6), and
$\iota$ the corresponding involution. Let $\EE \rightarrow \BB$ be
any algebraic family, with connected $\BB$ of course, such that
$\MM$ is one of the fibres. By looking at the canonical linear
systems, one can show that the fibres of $\EE$ over a Zariski-open
subset of $\BB$ are also double covers of smooth quadrics. Suppose
that we take the line bundle $\LL \rightarrow \BB$ which is some
high power of the fibrewise canonical bundle; it then follows that
the image of \eqref{eq:bb} consists of elements which commute with
$[\iota]$. Donaldson asked whether {\em all} symplectic
automorphisms act on $H_*(\MM)$ in a $\iota$-equivariant way; this
remains an open question.
\end{example}

As this concrete example suggests, there are currently no tools
strong enough to compute symplectic mapping class groups for
algebraic surfaces (or general symplectic four-manifolds) which are
not rational or ruled. However, the relation between $\pi_1(\BB)$ and
$\pi_0(Aut(\MM))$ can be probed by looking at the behaviour of some
particularly simple classes of symplectic automorphisms, and one of
these will be the subject of these lectures.

Namely, let $M$ be a closed symplectic four-manifold, and $L \subset
M$ an embedded Lagrangian two-sphere. One can associate to this a
Dehn twist or Picard-Lefschetz transformation, which is an element
$\tau_L \in Aut(M)$ determined up to isotopy. The definition is a
straightforward generalization of the classical Dehn twists in two
dimensions. However, the topology turns out to be rather different:
because of the Picard-Lefschetz formula
\begin{equation} \label{eq:p-l}
 (\tau_L)_*(x) = \begin{cases} x + (x \cdot l)\, l & x \in H_2(M;\Z), \\
 x & x \in H_k(M;\Z),\; k \neq 2 \end{cases}
\end{equation}
where $l = \pm [L]$ satisfies $l \cdot l = -2$, the square
$\tau_L^2$ acts trivially on homology, and in fact it is isotopic
to the identity in $\Diff(M)$. The obvious question is whether the
same holds in $Aut(M)$ as well. The first case which comes to mind
is that of the anti-diagonal in $M = S^2 \times S^2$ with the
monotone symplectic structure, and there $\tau_L^2$ is indeed
symplectically isotopic to the identity. But this is a rather
untypical situation: we will show that under fairly weak
conditions on a symplectic four-manifold, $[\tau_L^2] \in
\pi_0(Aut(M))$ is nontrivial whatever the choice of $L$. To take a
popular class of examples,

\begin{theorem} \label{th:one}
Let $M \subset \CP{n+2}$ be a smooth complete intersection of complex
dimension two. Suppose that $M$ is neither $\CP{2}$ nor $\CP{1}
\times \CP{1}$, which excludes the multidegrees $(1,\dots,1)$ and
$(2,1,\dots,1)$. Then the homomorphism
\begin{equation} \label{eq:inclusion}
\pi_0(Aut(M)) \longrightarrow \pi_0(\Diff(M))
\end{equation}
induced by inclusion is not injective.
\end{theorem}

There is in fact a slightly subtler phenomenon going on, which has to
do with the change in topology of $Aut(M)$ as the symplectic
structure varies. Let $\phi$ be a symplectic automorphism with
respect to the given symplectic form $\o$. We say that $\phi$ is {\em
potentially fragile} if there is a smooth family $\o^s$ of symplectic
forms, $s \in [0;\epsilon)$ for some $\epsilon>0$, and a smooth
family $\phi^s$ of diffeomorphisms such that $(\phi^s)^*\o^s = \o^s$,
with the following properties: (1) $(\phi^0,\o^0) = (\phi,\o)$; (2)
for all $s>0$, $\phi^s$ is isotopic to the identity inside
$Aut(M,\o^s)$. If in addition, (3) $\phi$ is not isotopic to the
identity in $Aut(M,\o)$, we say that $\phi$ is {\em fragile}. It is a
basic fact that squares of Dehn twists are always potentially
fragile, and so we have:

\begin{cor} \label{th:two}
Every two-dimensional complete intersection other than $\CP{2}$,
$\CP{1} \times \CP{1}$ admits a fragile symplectic automorphism.
\end{cor}

As suggested by their alternative name, Dehn twists do occur as
monodromy maps in families of algebraic surfaces, so the
nontriviality of $\tau^2$ proves that symplectic mapping class groups
do detect certain kinds of elements of $\pi_1$ of a moduli space,
which are hidden from ordinary topology. Moreover, the fragility
phenomenon has a natural interpretation in these terms.

Theorem \ref{th:one} and Corollary \ref{th:two} are taken from the
author's Ph.D.\ thesis \cite{seidel97}. Time and \cite{seidel01} have
made many of the technical arguments standard, and that frees us to
put more emphasis on examples and motivation, but otherwise the
structure and limitations of the original exposition have been
preserved. However, it seems reasonable to point out some related
results that have been obtained since then. For $K3$ and Enriques
surfaces containing two disjoint Lagrangian spheres $L_1,L_2$, it was
shown in \cite{seidel99} that $[\tau_{L_1}] \in \pi_0(Aut(M))$ has
infinite order, and therefore that the map \eqref{eq:inclusion} has
infinite kernel. \cite{khovanov-seidel98} proves that for the
noncompact four-manifold $M$ given by the equation $xy + z^{m+1} = 1$
in $\C^3$, there is a commutative diagram
\[
\xymatrix{
 B_m \ar[rr] \ar[d] && \pi_0(Aut^c(M)) \ar[d] \\
 S_m \ar[rr] && \pi_0(\Diff^c(M))
}
\]
where the upper $\rightarrow$ is injective. In particular, the
kernel of the right $\downarrow$ contains a copy of the pure braid
group $PB_m$. Similar phenomena happen for closed four-manifolds:
for instance, for a suitable symplectic form on the $K3$ surface,
one can show using \cite{khovanov-seidel98} and
\cite{seidel-thomas99} that the kernel of \eqref{eq:inclusion}
contains a copy of $PB_m$ for at least $m = 15$. In fact, a
simplified version of the same phenomenon (with a more direct
proof) already occurs for the del Pezzo surface $\CP{2} \# 5
\overline{\CP{2}}$, see Example \ref{th:dp5} below. All of this
fits in well with the idea that maps \eqref{eq:bb} should be an
important ingredient in understanding symplectic mapping class
groups of algebraic surfaces.
%
%

\acknowledgments Obviously, this work owes enormously to Simon
Donaldson, whose student I was while carrying it out originally.
Peter Kronheimer pointed out simultaneous resolution, which was a
turning-point. I have profited from discussions with Norbert A'Campo,
Denis Auroux, Mike Callahan, Allen Knutson, Anatoly Libgober, Dusa
McDuff, Leonid Polterovich, Ivan Smith, Dietmar Salamon, Richard
Thomas, Andrey Todorov and Claude Viterbo. Finally, I'd like to thank
the organizers and audience of the C.I.M.E.\ summer school for giving
me the opportunity to return to the subject.

\section{Definition and first properties}

\subsection{}
The construction of four-dimensional Dehn twists is standard
\cite{arnold95,seidel98b,seidel01}, but we will need the details as a
basis for further discussion. Consider $T^*\!S^2$ with its standard
symplectic form $\o$, in coordinates
\[
 T^*\!S^2 = \{ (u,v) \in \R^3 \times \R^3 \suchthat \leftsc u,v
 \rightsc = 0, ||v|| = 1\}, \quad
 \o = du \wedge dv.
\]
This carries the $O(3)$-action induced from that on $S^2$. Maybe less
obviously, the function $h(u,v) = ||u||$ induces a Hamiltonian circle
action $\sigma$ on $T^*\!S^2 \setminus S^2$,
\[
 \sigma_t(u,v) = \Big(
 \cos(t)u - \sin(t) ||u|| v,
 \cos(t)v + \sin(t) \frac{u}{||u||} \Big).
\]
$\sigma_\pi$ is the antipodal map $A(u,v) = (-u,-v)$, while for $t
\in (0;\pi)$, $\sigma_t$ does not extend continuously over the
zero-section. Geometrically with respect to the round metric on
$S^2$, $\sigma$ is the {\em normalized geodesic flow}, transporting
each tangent vector at unit speed (irrespective of its length) along
the geodesic emanating from it. Thus, the existence of $\sigma$ is
based on the fact that all geodesics on $S^2$ are closed. Now take a
function $r: \R \rightarrow \R$ satisfying $r(t) = 0$ for $t \gg 0$
and $r(-t) = r(t) - t$. The Hamiltonian flow of $H = r(h)$ is
$\phi_t(u,v) = \sigma_{t\,r'(||u||)}(u,v)$, and since $r'(0) = 1/2$,
the time $2\pi$ map can be extended continuously over the
zero-section as the antipodal map. The resulting compactly supported
symplectic automorphism of $T^*\!S^2$,
\[
 \tau(u,v) = \begin{cases}
 \sigma_{2\pi\,r'(||u||)}(u,v) & u \neq 0, \\
 (0,-v) & u = 0
 \end{cases}
\]
is called a model Dehn twist. To implant this local model into a
given geometric situation, suppose that $L \subset M$ is a
Lagrangian sphere in a closed symplectic four-manifold, and choose
an identification $i_0: S^2 \rightarrow L$. The Lagrangian tubular
neighbourhood theorem \cite[Theorem 1.5]{auroux-smith} tell us
that $i_0$ extends to a symplectic embedding \[i:
T^*_{\scriptscriptstyle \leq \lambda}S^2 \longrightarrow L\] of
the space $T^*_{\scriptscriptstyle \leq \lambda} S^2 \subset
T^*\!S^2$ of cotangent vectors of length $\leq \lambda$, for some
small $\lambda>0$. By choosing $r(t) = 0$ for $t \geq \lambda/2$,
one gets a model Dehn twist $\tau$ supported inside that subspace,
and then one defines the Dehn twist $\tau_L$ to be
\[
 \tau_L(x) = \begin{cases} i \tau i^{-1}(x) & x \in im(i), \\
 x & \text{otherwise.}
 \end{cases}
\]
The construction is not strictly unique, but it is unique up to
symplectic isotopy. The only choice that carries any topology is
the identification $i_0$, but this can be dealt with by observing
that $\tau$ is $O(3)$-equivariant, and $\Diff(S^2) \htp O(3)$ by
Smale's theorem. In particular, $\tau_L$ does not depend on a
choice of orientation of $L$.

If the circle action $\sigma$ extended smoothly over the
zero-section, then we could write down a compactly supported
symplectic isotopy between $\tau^2$ and the identity by moving along
the orbits,
\begin{equation} \label{eq:isotopy}
 \psi_t(u,v) = \sigma_{4\pi t\, r'(||u||)}(u,v).
\end{equation}
This may seem a pointless remark, since $\sigma$ does not extend over
$S^2$, but it comes into its own after a perturbation of the
symplectic structure. Take the standard symplectic form on $S^2$,
$\beta_v(X,Y) = \leftsc v, X \times Y \rightsc$, and pull it back to
$T^*\!S^2$. Then $\o^s = \o + s\beta$, $s \in \R$, is still an
$SO(3)$-invariant symplectic form.

\begin{prop} \label{th:local-fragility}
There is a smooth family $(\phi^s)$ of compactly supported
diffeomorphisms of $T^*\!S^2$, with the following properties: (1)
$\phi^s$ is symplectic for $\o^s$; (2) for all $s \neq 0$, $\phi^s$
is isotopic to the identity by an isotopy in $Aut^c(T^*\!S^2,\o^s)$;
(3) $\phi^0$ is the square $\tau^2$ of a model Dehn twist.
\end{prop}

We begin with an elementary general fact. For concreteness, we will
identify $\mathfrak{so}_3^* \iso \mathfrak{so}_3 \iso \R^3$ by using
the cross-product and the standard invariant pairing.

\begin{lemma} \label{th:so3}
Let $M$ be a symplectic manifold, carrying a Hamiltonian
$SO(3)$-action $\rho$ with moment map $\mu$. Then $h = ||\mu||$ is
the Hamiltonian of a circle action on $M \setminus \mu^{-1}(0)$.
\end{lemma}

\proof $h$ Poisson-commutes with all components of $\mu$ (since this
is true for the Poisson bracket on $\mathfrak{so}_3^*$, a well-known
fact from mechanics), so its flow maps each level set $\mu^{-1}(w)$
to itself. The associated vector field $X$ satisfies
\[
 X|\mu^{-1}(w) = K_{w/||w||}|\mu^{-1}(w)
\]
where $K$ are the Killing vector fields, which is clearly a circle
action (the quotient $\mu^{-1}(w)/S^1$ can be identified with the
symplectic quotient $M /\!/ SO(3)$ with respect to the coadjoint
orbit of $w$). \qed

The moment map for the $SO(3)$-action on $T^*\!S^2$ is $\mu(u,v) = -u
\times v$, so the induced circle action is just $\sigma$. With
respect to the deformed symplectic structures $\o^s$, the
$SO(3)$-action remains Hamiltonian but the moment map is $\mu^s(u,v)
= - sv - u \times v$, which is nowhere zero and hence gives rise to a
circle action $\sigma^s$ on the whole cotangent space. As $r
\rightarrow 0$, $\sigma^s$ converges on compact subsets of $T^*\!S^2
\setminus S^2$ to $\sigma$. For simplicity, assume that our model
Dehn twist $\tau$ is defined using a function $h$ which satisfies
$h'(t) = 1/2$ for small $t$. Then
\[
\phi^s(u,v) = \sigma^s_{4 \pi h'(||v||)}(u,v)
\]
for $s \neq 0$ defines a family of compactly supported
$\o_s$-symplectic automorphisms. These are all equal to the identity
in a neighbourhood of the zero section, hence they match up smoothly
with $\phi^0 = \tau_L^2$. By replacing $\sigma$ with $\sigma^s$ in
\eqref{eq:isotopy} one finds $\o^s$-symplectic isotopies between each
$\phi^s$, $s \neq 0$, and the identity. This concludes the proof of
Proposition \ref{th:local-fragility}. It is no problem to graft this
local construction into any Dehn twist, which yields:

\begin{cor} \label{th:fragility}
For any Lagrangian sphere $L$ in a closed symplectic four-manifold
$M$, the square $\tau_L^2$ of the Dehn twist is potentially fragile.
\qed
\end{cor}

\subsection{}
It is easy to see that any compactly supported $O(3)$-equivariant
symplectic automorphism of $T^*\!S^2$ has the form
\[
 \phi(x) = \begin{cases}
 \sigma_{2\pi r'(||x||)}(x) & x \notin S^2, \\
 A^k(x) & x \in S^2
 \end{cases}
\]
for some $k \in \Z$, and where $r: \R \rightarrow \R$ is a
function with $r(t) = 0$ for $r \gg 0$, and $r(-t) = r(t) - kt$
everywhere. There is no topologically nontrivial information in
this data except for $k$, so the space of such automorphisms is
weakly homotopy equivalent to the discrete set $\Z$ (the adjective
``weakly'' is a technical precaution, since there are several
slightly different choices for the topology on compactly supported
diffeomorphism groups, which however all have the same spaces of
continuous maps from finite-dimensional manifolds). We will now
see that the topology does not change if the equivariance
condition is dropped:

\begin{prop} \label{th:aut-cotangent}
The compactly supported automorphism group $Aut^c(T^*\!S^2)$ is
weakly homotopy equivalent to the discrete set $\Z$, with $1 \in \Z$
mapped to the model Dehn twist.
\end{prop}

In particular $[\tau^k] \in \pi_0(Aut^c(T^*\!S^2))$ is nontrivial for
all $k \neq 0$. The result also says that up to isotopy and
iterating, a Dehn twist is the only construction of a symplectic
automorphism that can be done locally near a Lagrangian sphere.

\proof This is an easy consequence of Gromov's work. Take $M = S^2
\times S^2$ with the standard product symplectic form (in which both
factors have the same volume), $L = \{x_1 + x_2 = 0\}$ the
antidiagonal, and $\Delta = \{x_1 = x_2\}$ the diagonal. Consider the
groups
\begin{align*}
 \G_1 &= \{ \phi \in Aut(M) \suchthat \phi(\Delta) = \Delta \}, \\
 \G_2 &= \{ \phi \in \G_1 \suchthat \phi|\Delta = id\}, \\
 \G_3 &= \{ \phi \in \G_2 \suchthat \phi|U = id \text{ for some
 open $U \supset \Delta$} \}.
\end{align*}
First of all, $M \setminus \Delta$ is isomorphic to
$T^*_{\scriptscriptstyle <\lambda}\!S^2$ some $\lambda$, with $L$
corresponding to the zero section. Therefore we have a weak
homotopy equivalence $\G_3 \htp Aut^c(T^*\!S^2)$. Next, there is a
weak fibration
\[
 \G_3 \longrightarrow \G_2 \stackrel{D}{\longrightarrow} Map(S^2,S^1)
\]
where $Map(S^2,S^1)$ is thought of as the group of unitary gauge
transformations of the normal bundle to $\Delta$, and $D$
essentially the map which associates to each automorphism its
derivative in normal direction. It is an easy observation that
$Map(S^2,S^1) \htp S^1$. Third, we have a weak fibration
\[
 \G_2 \longrightarrow \G_1 \longrightarrow \Diff^+(S^2),
\]
with $\Diff^+(S^2) \htp SO(3)$. Finally
\[
 \G_1 \longrightarrow \G_0 \longrightarrow {\mathcal S}_\Delta,
\]
where ${\mathcal S}_\Delta$ is the space of embedded symplectic
two-spheres in $S^2 \times S^2$ which can be mapped to $\Delta$ by a
symplectic automorphism. Gromov's theorem says that $Aut(M) \htp
(SO(3) \times SO(3)) \rtimes \Z/2$, and a variation of another of his
basic results is that ${\mathcal S}_\Delta \htp SO(3)$. Appying these
sequences in the reverse order, one finds that $\G_1$ is homotopy
equivalent to $SO(3) \times \Z/2$, and that $\G_2 \htp \Z/2$, so the
higher homotopy groups of $\G_3$ vanish while $\pi_0(\G_3)$ sits in a
short exact sequence
\begin{equation} \label{eq:z2z}
 1 \rightarrow \Z \longrightarrow \pi_0(\G_3) \stackrel{\alpha}{\longrightarrow}
 \Z/2 \rightarrow 1,
\end{equation}
where $\alpha$ assigns to each symplectomorphism the sign $\phi_*[L]
= \pm [L]$. The last step yields the following additional
information: take a map $\phi \in \G_3$ which preserves the
orientation of $L$, and let $\phi_t$ be a homotopy from it to the
identity inside $\G_2$. Then the element of $ker(\alpha)$ represented
by $\phi$ is the degree of $S^1 \rightarrow Sp_4(\R)$, $t \mapsto
D_x\phi_t$ at any point $x \in \Delta$. By applying this to the
isotopy $\tau_L^2 \htp id$ constructed in Example \ref{th:s2s2}
below, one sees that \eqref{eq:z2z} does not split, and that
$[\tau_L]$ is a generator of $\pi_0(Aut^c(T^*\!S^2)) \iso \Z$. \qed

It is an interesting exercise to see how the above argument changes
if one passes to the symplectic form $\o^s = \o + s\beta$ for small
$s$.

\subsection{}
Corollary \ref{th:fragility} is too essential to pass it off as the
result of some ad hoc local construction. A proper understanding
involves looking at the real nature of Dehn twists as monodromy maps.

\begin{defn} \label{th:lefschetz}
Let $S$ be an oriented surface, possibly non-compact or with
boundary. A (six-dimensional) symplectic Lefschetz fibration over
$S$ is a six-manifold $E$ with a proper map $\pi: E \rightarrow
S$, $\pi^{-1}(\partial S) = \partial E$, a closed two-form $\Omega
\in \Omega^2(E)$, a complex structure $J_E$ defined on a
neighbourhood of the set of critical points $E^{crit}$, and a
positively oriented complex structure $j_S$ defined on a
neighbourhood of the set of critical values $S^{crit}$. The
requirements are:
\begin{itemize} \itemsep1em
\item
Near the critical points, $\pi$ is a holomorphic map with respect to
$J_E$ and $j_S$, and the critical points themselves are
nondegenerate. Moreover, $E^{crit}$ is disjoint from $\partial E$,
and $\pi|E^{crit}$ is injective.
\item
$\Omega$ is a K{\"a}hler form for $J_E$ in a neighbourhood of
$E^{crit}$. For any point $x \notin E^{crit}$, the restriction of
$\Omega_x$ to $TE_x^v = \ker(D\pi_x)$ is nondegenerate.
\end{itemize}
\end{defn}

The geometry of these fibrations is not very different from the
familiar four-dimensional case treated in \cite{auroux-smith}; one
possible reference for the results stated below is
\cite{seidel01}. Away from the critical fibres they are symplectic
fibrations, and in fact carry a preferred Hamiltonian connection
$TE^h$, the $\Omega$-orthogonal complement to $TE^v$ (the word
``Hamiltonian'' refers to the structure group $Aut^h$ from Remark
\ref{th:aut-h}, and does {\em not} mean that the monodromy
consists of maps Hamiltonian isotopic to the identity). Hence, for
any smooth path $\gamma: [0;1] \rightarrow S \setminus S^{crit}$
we have a canonical parallel transport map $P_\gamma:
E_{\gamma(0)} \rightarrow E_{\gamma(1)}$. Given a path $\gamma:
[0;1] \rightarrow S$ with $\gamma^{-1}(S^{crit}) = \{1\}$,
$\gamma'(1) \neq 0$, one can look at the limit of
$P_{\gamma|[0;t]}$ as $t \rightarrow 1$, and this gives rise to a
Lagrangian two-sphere $V_\gamma \subset E_{\gamma(0)}$, which is
the vanishing cycle of $\gamma$. The Picard-Lefschetz theorem says
that if $\lambda$ is a loop in $S \setminus S^{crit}$ with
$\lambda(0) = \lambda(1) = \gamma(0)$, winding around $\gamma$ in
positive sense, its monodromy is the Dehn twist around the
vanishing cycle, at least up to symplectic isotopy:
\[ P_\lambda \htp \tau_{V_\gamma} \in Aut(E_{\gamma(0)}). \]

Let's pass temporarily to algebro-geometric language, so $\pi: \EE
\rightarrow \SS$ is a proper holomorphic map from a threefold to a
curve, with the same kind of critical points as before, and $\LL
\rightarrow \EE$ is a relatively very ample line bundle. Atiyah
\cite{atiyah58} (later generalized by Brieskorn \cite{brieskorn68})
discovered the phenomenon of {\em simultaneous resolution}, which can
be formulated as follows: let $r: \hat{\SS} \rightarrow \SS$ be a
branched covering which has double ramification at each preimage of
points in $\SS^{crit}$. Then there is a commutative diagram
\[
\xymatrix{
 \hat{\EE} \ar[rr]^-{R} \ar[d]^{\hat{\pi}} && \EE \ar[d]^{\pi} \\
 \hat{\SS} \ar[rr]^-{r} && \SS
}
\]
where $\hat{\pi}$ has no critical points (proper smooth morphism),
and the restriction of $R$ gives an isomorphism \[\hat{\EE}
\setminus R^{-1}(\EE^{crit}) \stackrel{\iso}{\longrightarrow}
r^*(\EE \setminus \EE^{crit}).\] In particular, away from the
singular fibres $\hat{\EE}$ is just the pullback of $\EE$. If
$\lambda$ is a small loop in $\SS \setminus \SS^{crit}$ going once
around a critical value, then its iterate $\lambda^2$ can be
lifted to $\hat{\SS}$, which means that the monodromy around it
must be isotopic to the identity as a diffeomorphism. Of course,
by the Picard-Lefschetz formula $P_{\lambda^2} \htp \tau_V^2$ for
the appropriate vanishing cycle $V$. The preimage of each critical
point $x \in \EE$ is a rational curve $C_x \subset \hat{\EE}$ with
normal bundle $\OO(-2)$ in its fibre. Suppose that there is a line
bundle $\Lambda \rightarrow \hat{\EE}$ such that $\Lambda|C_x$ has
positive degree for each $x$ (this may or may not exist, depending
on the choice of resolution). Then $\hat{\LL} = \LL^{\otimes d}
\otimes \Lambda^{\otimes e}$ is relatively very ample for $d \gg e
\gg 0$. This shows that the monodromy around $\lambda^2$ becomes
symplectically trivial after a change of the symplectic form,
which is essentially the same property as potential fragility of
$\tau_V^2$ except that algebraic geometry does not actually allow
us to see this change as a continuous deformation. However, one
can easily copy the local construction of the simultaneous
resolution in the symplectic setting, and this gives an
alternative proof of Corollary \ref{th:fragility} avoiding any
explicit computation.

\begin{remark}
More generally, potential fragility occurs naturally in situations
involving hyperk{\"a}hler quotients. Let $X$ be a hyperk{\"a}hler
manifold, and pick a preferred complex structure on it. Suppose that
it carries a hyperk{\"a}hler circle action with moment map $h =
(h_\R,h_\C): X \rightarrow \R \times \C$, and a connected component
of the fixed point set on which $h \equiv 0$. For simplicity we will
assume that the action is otherwise free, and ignore problems arising
from the noncompactness of $X$ (so the following statements are not
entirely rigorous). If one fixes $s \in \R$ then
\[
 X^s_{\C^*} = (h_\R^{-1}(s) \setminus h_\C^{-1}(0)) /S^1
 \stackrel{h_\C}{\longrightarrow} \C^*
\]
is a holomorphic map, and the total space carries a natural
quotient K{\"a}hler form. One can therefore define the monodromy
around a circle of some radius $\epsilon>0$ in the base $\C^*$,
which is a symplectic automorphism $\phi^s$ of the hyperk{\"a}hler
quotient $X^s_{\epsilon} = (h_\R^{-1}(s) \cap
h_{\C}^{-1}(\epsilon))/S^1$. Varying $s$ does not affect the
complex structure on $X^s_\epsilon$, but the K{\"a}hler class
varies, so one can consider the $\phi^s$ as a family of
automorphisms for a corresponding family $\o^s$ of symplectic
forms on a fixed manifold. For $s = 0$ there is a singular fibre
$X^0_0$ at the center of the circle, and one would hope that
$\phi^0$ reflects this fact; in contrast, for $s \neq 0$ we have
that $X^s_z$ is smooth for all $z$, so $\phi^s$ is symplectically
isotopic to the identity. The case of squared Dehn twists is a
particularly simple example of this, with $X = \mathbb{H}^2$; see
\cite{kronheimer89}. A straightforward generalization leads to
analogues of $\tau^2$ on $T^*\CP{n}$, which were discussed in
\cite{seidel99}.
\end{remark}

\subsection{}
Donaldson's theory of almost holomorphic functions is an attempt
to reduce all questions about symplectic four-manifolds to
two-dimensional ones, and hence to combinatorial group theory. The
paper \cite{donaldson02} achieves this for the fundamental
classification problem, but the wider program also embraces
symplectic mapping groups. The relevant deeper results are still
being elaborated, but the elementary side of the theory is
sufficient to understand the potential fragility of squared Dehn
twists. The following discussion is due to Donaldson (except
possibly for mistakes introduced by the author). Compared to the
exposition in \cite{auroux-smith}, to which the reader is referred
for the basic theory of Lefschetz pencils, we will just need to
exercise a little more care concerning the definition of
symplectic forms on the total spaces.

Let $S$ be a closed oriented surface, equipped with a symplectic
form $\eta$ and a finite set of marked points $\Sigma =
\{z_1,\dots,z_p\}$, which may be empty. We assume that the Euler
characteristic $\chi(S \setminus \Sigma) < 0$. Denote by
$Aut^h(S,\Sigma)$ the group of symplectic automorphisms of $S$
which are the identity in a neighbourhood of $\Sigma$, with the
Hamiltonian topology. For any simple closed curve $\gamma \subset
S \setminus \Sigma$, we have the (classical) Dehn twist
$t_\gamma$, which is an element of $Aut^h(S,\Sigma)$ unique up to
isotopy within that topological group (note that if
$\gamma,\gamma'$ are nonseparating curves which are isotopic to
each other, but not Hamiltonian isotopic, then $t_\gamma$ and
$t_{\gamma'}$ have different classes in $\pi_0(Aut^h)$). Choose a
small loop $\zeta_k$ around each $z_k$. Take a finite ordered
family $(\gamma_1,\dots,\gamma_m)$ of simple closed
non-contractible curves in $S \setminus \Sigma$, such that
\begin{equation} \label{eq:mcg}
 t_{\gamma_1}\dots t_{\gamma_m} \htp t_{\zeta_1} \dots
 t_{\zeta_p}
\end{equation}
in $Aut^h(S,\Sigma)$. From this one constructs a four-manifold $M$
together with a family $\o^s$ of closed forms, which are symplectic
for $s \gg 0$. For brevity, we will call this an {\em asymptotically
symplectic manifold}. The first step is take the (four-dimensional
topological) Lefschetz fibration $\tilde{M} \rightarrow S^2$ with
smooth fibre $S$ and vanishing cycles $\gamma_1,\dots,\gamma_m$.
Using a suitable Hamiltonian connection, one can define a closed
two-form $\tilde{\o}$ on $\tilde{M}$ whose restriction to each smooth
fibre is symplectic. The family $\tilde{\o}^s = \tilde{\o} + s
\beta$, where $\beta$ is the pullback of a positive volume form on
$S^2$, consists of symplectic forms for $s \gg 0$. Each base point
$z_k$ will give rise to a section, whose image is a symplectic sphere
with self-intersection $-1$. Blowing down these spheres completes the
construction of $(M,\{\o^s\})$. Of course there is some choice in the
details, but the outcome is unique up to asymptotic symplectic
isomorphism, which is the existence of a family of diffeomorphisms
$\{\phi^s\}$ which are symplectic for $s \gg 0$; and moreover, this
family is canonical up to asymptotically symplectic isotopy, which is
enough for our purpose. For later reference, we note the following
fact about the cohomology class of $\o^s$. The primitive part
$H^2(M;\R)^{prim}$, which is just the quotient of $H^2(M;\R)$ by the
Poincar{\'e} dual of the fibre $S \subset M$, can be described as the
middle cohomology group of a complex \cite{donaldson98}
\begin{equation} \label{eq:donaldson-seq}
 H^1(S \setminus \Sigma;\R) \stackrel{a}{\longrightarrow} \R^m
 \stackrel{a'}{\longrightarrow} H^1(S;\R)
\end{equation}
where $a$ is given by integrating over the $\gamma_k$, and $a'$
involves a certain dual set of vanishing cycles $\gamma_k'$. The
class of $\o^s$ in $H^2(M;\R)^{prim}$ is independent of $s$, and is
represented by a vector in $\R^m$ in \eqref{eq:donaldson-seq} defined
by choosing a one-form $\theta$ on $S \setminus \Sigma$ with $d\theta
= \eta$, and integrating that over the $\gamma_k$. In particular, if
$\theta$ can be chosen in such a way that $\int_{\gamma_k} \theta =
0$ for all $k$, then all $\o^s$ are multiples of $PD([S])$, which is
the case of a Lefschetz {\em pencil}.

If one replaces the $\gamma_k$ by curves Hamiltonian isotopic to
them, $M$ remains the same, up to the same kind of isomorphism as
before. We call the equivalence class of $(\gamma_1,\dots,\gamma_m)$
under this relation a {\em Lefschetz fibration datum}; this will be
denoted by $\Gamma$, and the associated manifold by
$(M_\Gamma,\{\o_{\Gamma}^s\})$. More interestingly, there are two
nontrivial modifications of a Lefschetz fibration datum which do not
change $M$; together they amount to an action of $G =
\pi_0(Aut^h(S,\Sigma)) \times B_m$ on the set of such data. The first
factor acts by applying a symplectic automorphism $\phi$ to all of
the $\gamma_k$, and the generators of the braid group $B_m$ act by
elementary Hurwitz moves
\begin{equation} \label{eq:hurwitz}
 (\gamma_1,\dots,\gamma_m) \longmapsto
 (\gamma_1,\dots,\gamma_{k-1},t_{\gamma_k}(\gamma_{k+1}),\gamma_k,
 \gamma_{k+2},\dots,\gamma_m).
\end{equation}
Roughly speaking, what the two components of the $G$-action do is to
change the way in which the fibre of $\tilde{M}_\Gamma$ is identified
with $S$, respectively the way in which its base is identified with
$S^2$. By uniqueness, we have for every $g \in G$ such that
$g(\Gamma) = \Gamma$ an induced asymptotically symplectic
automorphism $\{\phi^s\}$ of $M_\Gamma$. Denoting by $G_\Gamma
\subset G$ the subgroup which stabilizes $\Gamma$, and by
$Aut(M_\Gamma,\{\o_\Gamma^s\})$ the group of asymptotically
symplectic automorphisms, we therefore have a canonical map
\begin{equation} \label{eq:lefschetz-aut}
 G_\Gamma \longrightarrow \pi_0(Aut(M_\Gamma,\{\o_\Gamma^s\})).
\end{equation}
(in the case of a Lefschetz pencil, the right hand side reduces to
$Aut(M_\Gamma,\o_{\Gamma}^\sigma)$ for some fixed $\sigma \gg 0$).
Usually \eqref{eq:lefschetz-aut} is not injective. For instance,
consider the situation where two subsequent curves $\gamma_k$,
$\gamma_{k+1}$ are disjoint. Applying \eqref{eq:hurwitz} just
exchanges the curves; the square of this operation is a nontrivial
element of $G_\Gamma$, but the associated asymptotically symplectic
automorphism is isotopic to the identity. This can be most easily
seen by thinking of families of Lefschetz fibrations: in our case, we
have a family parametrized by $S^1$ in which two critical values in
$S^2$ rotate around each other, and whose monodromy is the image of
our Hurwitz move in \eqref{eq:lefschetz-aut}; but since the vanishing
cycles are disjoint, we can move the two critical points into the
same fibre, and so the family can be extended over $D^2$, which
trivializes the monodromy.

Suppose that we are in the Lefschetz pencil situation where
$\int_{\gamma_k} \theta = 0$, and that two subsequent curves
$\gamma_k$, $\gamma_{k+1}$ agree. One can then use their bounding
``Lefschetz thimbles'' to construct a Lagrangian sphere $L \subset
M$, and its inverse Dehn twist $\tau_L^{-1}$ is the image of the
elementary Hurwitz move \eqref{eq:hurwitz} under
\eqref{eq:lefschetz-aut}. Now move $\gamma_k$, $\gamma_{k+1}$ away
from each other in a non-Hamiltonian way, by an opposite amount of
area. The resulting new configuration of curves
$\gamma_1',\dots,\gamma_m'$ still satisfies the basic equation
\eqref{eq:mcg}, and defines the same four-manifold $M_{\Gamma'} =
M_\Gamma$ with a different symplectic form: an argument using
\eqref{eq:donaldson-seq} shows that $\o_{\Gamma'}^s$ differs from
$\o_{\Gamma}^s$ by a multiple of $PD(L)$, which becomes comparatively
small as $s \rightarrow \infty$. Since $\gamma_k', \gamma_{k+1}'$ are
disjoint, the element of $G_\Gamma$ which led to $\tau_L^2$ now
becomes an element of $G_\Gamma'$ inducing a trivial asymptotically
symplectic automorphism, which is the statement of potential
fragility in this framework.

\begin{remark}
We should briefly mention the expected deeper results concerning the
map \eqref{eq:lefschetz-aut} (these were first stated by Donaldson,
and their proof is the subject of ongoing work of
Auroux-Munoz-Presas). The main idea is that the image of
\eqref{eq:lefschetz-aut} for Lefschetz pencils should ultimately
exhaust the symplectic automorphism group as the degree of the pencil
goes to $\infty$. More precisely, given a symplectic manifold and
integral symplectic form $\o$, and an arbitrary symplectic
automorphism $\phi$, there should be a Lefschetz pencil whose fibres
lie in the class $k[\o]$ for $k \gg 0$, and an element of the
resulting $G_\Gamma$ which maps to $[\phi]$. There is also a list of
relations for the kernel of \eqref{eq:lefschetz-aut} which is
conjectured to be complete in a suitable $k \rightarrow \infty$
sense, but a rigorous formulation of that would be quite complicated
since it involves ``degree doubling''.
\end{remark}

\subsection{}
As usual, let $L \subset M$ be a Lagrangian sphere in a closed
symplectic four-manifold. Having considered the fragility of
$\tau_L^2$ from different points of view, we now turn to the main
question, which is whether it is isotopic to the identity in
$Aut(M)$. We know that this is a nontrivial question because the
answer for the corresponding local problem is negative, by
Proposition \ref{th:aut-cotangent}, and as mentioned in the
Introduction this answer carries over to the vast majority of
closed four-manifolds. For now, however, the discussion will start
from the opposite direction, as we try to accumulate examples
where $\tau_L^2$ {\em is} symplectically isotopic to the identity,
and then probe the line where something nontrivial happens.

First of all, there is an elementary construction based directly
on the circle action $\sigma$ used in the definition of the Dehn
twist.

\begin{lemma} \label{th:extend-action}
Suppose that there is a Hamiltonian circle action $\bar\sigma$ on $M
\setminus L$ and a Lagrangian tubular neighbourhood $i:
T^*_{\scriptscriptstyle <\lambda}S^2 \rightarrow M$ of $L$ which is
equivariant with respect to $\sigma$, $\bar\sigma$. Then $\tau_L^2$
is isotopic to the identity in $Aut(M)$. \qed
\end{lemma}

The proof is straightforward, and we leave it to the reader.

\begin{example} \label{th:s2s2}
As in the proof of Proposition \ref{th:aut-cotangent} take $M = S^2
\times S^2$ with the monotone symplectic form, and $L = \{x_1 + x_2 =
0\}$ the antidiagonal. The diagonal $SO(3)$-action has moment map
$\mu(x) = -x_1-x_2 \in \R^3$, and from Lemma \ref{th:so3} above we
know that $\bar{h}(x) = ||x_1+x_2||$ is the moment map for a circle
action $\bar\sigma$ on $M \setminus L$. This has the desired property
with respect to any $SO(3)$-equivariant Lagrangian tubular
neighbourhood for $L$. A slight refinement of Lemma
\ref{th:extend-action} shows that $\tau_L$ itself is symplectically
isotopic to the involution $(x_1,x_2) \mapsto (x_2,x_1)$. Somewhat
less transparently, this could also be derived from Gromov's Theorem
\ref{th:gromov}.
\end{example}

\begin{example} \label{th:linkages}
A related case is the ``regular pentagon space'', a manifold often
used as a basic example in the theory of symplectic quotients
\cite[Chapter 4 \S 5]{newstead} \cite[Chapter 16.1]{kirwan84}
\cite[Chapter 8]{mumford-fogarty-kirwan} (incidentally, it is also
the same as the Deligne-Mumford space $\overline{\mathcal M}_{0,5}$).
Take $S^2$ with its standard symplectic form, and consider the
diagonal action of $SO(3)$ on $(S^2)^5$ with moment map $\mu(x) =
-(x_1 + \dots + x_5)$. The symplectic quotient
$
 M = \mu^{-1}(0)/SO(3)
$ is the space of quintuples of vectors of unit length in $\R^3$
which add up to zero, up to simultaneous rotation. This is a
compact symplectic four-manifold, and it contains a natural
Lagrangian sphere
\begin{equation} \label{eq:l1-sphere}
 L_1 = \{x_1 + x_2 = 0\}.
\end{equation}
$M \setminus L_1$ carries a Hamiltonian circle action $\bar\sigma_1$,
given by rotating $x_1$ around the axis formed by $x_1 + x_2$ while
leaving $x_1+x_2$, $x_3$, $x_4$, $x_5$ fixed. The relevant moment map
is $\bar{h}_1(x) = ||x_1+x_2||$ as before, which already looks much
like our standard circle action on $T^*\!S^2 \setminus S^2$. Indeed,
one can find a tubular neighbourhood of $L_1$ satisfying the
conditions of Lemma \ref{th:extend-action}, so $\tau_{L_1}^2$ is
symplectically isotopic to the identity. In fact, by cyclically
permuting coordinates, one finds a configuration of Lagrangian
spheres $L_1,\dots,L_5$ whose intersections are indicated by a
pentagon graph
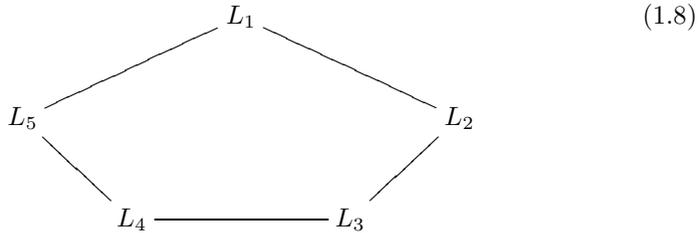
\begin{equation} \label{eq:a5-cyclic}
\xymatrix{
 && L_1 \ar@{-}[drr] \ar@{-}[dll] && \\
 L_5 \ar@{-}[dr] &&&& L_2 \ar@{-}[dl] \\
 & L_4 \ar@{-}[rr] && L_3 &
}
\end{equation}
Because of the resulting braid relations \cite[Appendix]{seidel98b},
$\tau_{L_1},\dots,\tau_{L_4}$ generate a homomorphism $B_5
\rightarrow \pi_0(Aut(M))$; on the other hand, we have the additional
relation $[\tau_{L_k}^2] = 1$, so this actually factors through the
symmetric group $S_5$.

It is worth while to identify $M$ more explicitly. Take the maps
induced by inclusion $j$ and projection $p$,
\[
H^2((S^2)^5;\R) \stackrel{j^*}{\longrightarrow} H^2(\mu^{-1}(0);\R)
\stackrel{p^*}{\longleftarrow} H^2(M;\R).
\]
Our group being $SO(3)$, a look at the standard spectral sequence
shows that cohomology and equivariant cohomology coincide in degree
two. This implies that $p^*$ is an isomorphism. Now, the pullback of
the symplectic form on $M$ via $p$ agrees with the restriction of the
symplectic form on $(S^2)^5$ via $j$, and the same holds for the
first Chern classes of their respective tangent bundles. We conclude
that $M$ is monotone, so by general classification results
\cite{lalonde-mcduff95d} it must be either $\CP{1} \times \CP{1}$ or
$\CP{2}$ blown up at $0 \leq k \leq 8$ points. The same consideration
with equivariant cohomology as before, together with Kirwan's
surjectivity theorem, shows that $j^*$ is onto, so $b_2(M) \leq 5$.
On the other hand, by looking at the intersection matrix of the
configuration \eqref{eq:a5-cyclic} one sees that the part of
$H^2(M;\R)$ orthogonal to the symplectic class has at least dimension
4. Therefore $b_2(M) = 5$ and so
\[
M \iso \CP{2} \# 4 \overline{\CP{2}}.
\]
A more elementary approach is to observe that $M \setminus (L_1 \cup
L_3)$ carries a $T^2$-action with three fixed points, which directly
yields $\chi(M) = 7$. Finally, one can vary this example by
considering quintuples of vectors of different lengths, see
\cite{klyachko94,hausmann-knutson97, hausmann-knutson98}. This yields
examples of Lagrangian spheres on $\CP{2} \# 2\overline{\CP{2}}$ and
$\CP{2} \# 3\overline{\CP{2}}$ with $\tau^2$ symplectically isotopic
to the identity, however the relevant symplectic forms are not
monotone.
\end{example}

Another way of finding examples of $\tau^2 \htp id$ is based on
the connection with algebraic geometry and monodromy, which means
on the construction of suitable families of algebraic surfaces
together with the Picard-Lefschetz theorem.

\begin{lemma} \label{th:no-monodromy}
Let $\pi: \EE \rightarrow \BB$ be a proper smooth map between
quasi-projective varieties of relative dimension 2, and $\LL
\rightarrow \EE$ a relatively very ample line bundle. Suppose that
there is a partial compactification $\bar{\pi}: \overline{\EE}
\rightarrow \overline{\BB}$ where the total space and base are
still smooth, and with the following properties: (1) The
discriminant $\Delta = \overline{\BB} \setminus \BB$ is a
hypersurface, and the fibre over a generic point $\delta \in
\Delta$ is reduced and has a single ordinary double point
singularity; (2) the meridian around $\delta$ is an element of
order 2 in $\pi_1(\BB)$; (3) $\LL$ extends to a relatively very
ample line bundle $\bar{\LL} \rightarrow \bar{\EE}$. If these
conditions hold, the smooth fibre $\EE_b$, $b \in \BB$, with its
induced K{\"a}hler structure, contains a Lagrangian sphere $L$
such that $\tau_L^2$ is symplectically isotopic to the identity.
\end{lemma}

\proof Take a generic point $\delta \in \Delta$, a neighbourhood
$U \subset \bar{\BB}$ of $\delta$, and local holomorphic
coordinates $(\zeta,y_1,\dots,y_n): U \rightarrow \C^{n+1}$ such
that $\Delta = \{\zeta = 0\}$. Take a small generic value of the
map $y \circ \bar{\pi}: \bar\EE|U \rightarrow U \rightarrow \C^n$.
The preimage of that is a smooth threefold $\bar\EE_y$ with a
holomorphic map $\bar\pi_y = \zeta \circ \bar\pi: \bar\EE_y
\rightarrow U_y \subset \C$, such that $\bar\pi_y^{-1}(0)$ has a
single ordinary double point. This implies that $\bar\pi_y$ has a
nondegenerate critical point. After using $\bar{\LL}$ to put a
suitable K{\"a}hler form on $\bar\EE_y$, we find that the
monodromy of $\bar\pi_y$ around $0$ is the Dehn twist along some
Lagrangian sphere. On the other hand, the monodromy is the image
of the meridian around $\delta$ by the map \eqref{eq:bb}, so it
must be of order two. \qed

\begin{example} \label{th:dp1}
A classical case is that of del Pezzo surfaces with small rank. The
necessary algebro-geometric background can be found in the first few
lectures of \cite{demazure-pinkham-teissier77}. Fix $2 \leq k \leq
4$. Consider a configuration $b = \{b_1,\dots,b_k\}$ of $k$ unordered
distinct points in $\CP{2}$ which are in general position, meaning
that no three of them are collinear. If that is the case, the
anticanonical bundle on the blowup $Bl_b(\CP{2})$ is ample, in fact
very ample. Over the space $\BB \subset Conf_k(\CP{2})$ of
configurations $c$ in general position, there is a natural family of
blowups $\EE \rightarrow \BB$, and the anticanonical bundle is
relatively very ample. The action of $PSL_3(\C)$ on $\BB$ is
transitive, and therefore $\pi_1(\BB)$ is the $\pi_0$ of the
stabilizer of any point, in fact
\begin{equation} \label{eq:small-pi1}
\pi_1(\BB) \iso S_k.
\end{equation}
What this says is that a symplectic automorphism obtained as
monodromy map from the family of blowups is symplectically isotopic
to the identity iff it acts trivially on homology.

To partially compactify $\BB$, we will now relax the genericity
conditions by allowing two points to collide. Let $b \subset
\OO_{\CP{2}}$ be an ideal sheaf of length $k$. This means that it is
a configuration of points with multiplicities, which add up to $k$,
and additional infinitesimal information at the multiple points. We
say that $b$ is in almost general position if any point occurs at
most with multiplicity two, and its restriction to any line in
$\CP{2}$ has length at most two. The space of such ideals is a
partial compactification $\bar{\BB}$ of $\BB$, and the discriminant
$\Delta$ is a smooth divisor. For $\delta \in \Delta$,
$Bl_\delta(\CP{2})$ is a surface with an ordinary double point. This
means that $\EE$ extends to a family $\bar{\EE} \rightarrow
\bar{\BB}$ such that the fibres over $\Delta$ have an ordinary double
point. It is not difficult to show that $\bar{\EE}$ is smooth, and
that the anticanonical bundle is still fibrewise very ample. In view
of \eqref{eq:small-pi1}, Lemma \ref{th:no-monodromy} implies that
there is a Lagrangian sphere $L$ in the smooth fibre, which is
$\CP{2} \# k \overline{\CP{2}}$ with a monotone symplectic form, such
that $\tau_L^2$ is symplectically isotopic to the identity (of
course, for $k = 4$ we already know this from Example
\ref{th:linkages}, but the cases $k = 2,3$ are new).
\end{example}

The structure of the moduli space changes for del Pezzo surfaces
of rank $5 \leq k \leq 8$, where the action of $PSL_3(\C)$ on the
corresponding space of generic configurations is no longer
transitive. As one would expect from the general philosophy, this
also affects the structure of symplectic mapping class groups. The
first case $k = 5$ can be treated by elementary means, and we will
do so now.

\begin{example} \label{th:dp5}
Take $\BB^{ord} \subset Conf_5^{ord}(\CP{2})$ to be the space of
ordered quintuples of points $b = (b_1,\dots,b_5)$ in the
projective plane which are in general position, in the same sense
as before. $PSL_3(\C)$ acts freely on this, and the quotient
$\BB^{ord}/PSL_3(\C)$ is isomorphic to the moduli space of ordered
quintuples of points on the line,
$Conf_5^{ord}(\CP{1})/PSL_2(\C)$. One can see this by direct
computation: each $PSL_3(\C)$-orbit on $\BB^{ord}$ contains
exactly one point of the form
\begin{equation} \label{eq:b-generic}
\begin{aligned}
 & b_1 = [1:0:0], \; b_2 = [0:1:0], \;
 b_3 = [0;0;1], \; \\ & \qquad b_4 = [1:1:1], \;
 b_5 = [z/w;(1-z)/(1-w);1]
\end{aligned}
\end{equation}
with $z,w \neq \{0,1\}$ and $z \neq w$; and correspondingly, each
$PSL_2(\C)$-orbit on $\CP{1}$ contains a unique configuration of
the form $(0,1,\infty,z,w)$. A more geometric construction goes as
follows: there is a unique (necessarily nonsingular, by the
general position condition) conic $Q$ which goes through the
points $b_1,\dots,b_5$. One can identify $Q \iso \CP{1}$ and then
the $b_k \in Q$ become a configuration of points on the line. A
straightforward computation shows that this gives back our
previous identification; the requirement that this should work out
explains the strange coordinates used in \eqref{eq:b-generic}.

This approach can be imitated on a symplectic level. Let $M =
\CP{2} \# 5 \overline{\CP{2}}$ with its monotone symplectic
structure. Let $E_1,\dots,E_5$ be the homology classes of the
exceptional curves in the $\overline{\CP{2}}$ summands, and $L$
the homology class of the line in $\CP{2}$. Take an arbitrary (not
generic in any sense) compatible almost complex structure $J$.
Each of the classes
\begin{equation} \label{eq:e-classes}
E_1,\dots,E_5,2L-E_1-\cdots-E_5
\end{equation}
is minimal, in the sense that it cannot be written as the sum of
two classes of positive symplectic area, hence the moduli space of
$J$-holomorphic spheres in that class is compact. The adjunction
formula \cite[Theorem 2.2.1]{mcduff-pos} proves that this space
consists of embedded spheres, and the regularity theorem from
\cite{hofer-lizan-sikorav} implies that it is smooth. By deforming
to the standard complex structure, one sees that each class
\eqref{eq:e-classes} is represented by a unique embedded
$J$-holomorphic sphere. The multiplicity theorem \cite[Theorem
2.1.1]{mcduff-pos} shows that the sphere representing $2L - E_1 -
\dots - E_5$ intersects each $E_k$ sphere transversally in a
single point. Hence, by identifying that sphere with $\CP{1}$ one
gets an element of $Conf_5^{ord}(\CP{1})/PSL_2(\C)$. This can be
done on the fibres of a symplectic fibration, as long as the
homological monodromy is trivial (allowing one to identify the
homology classes \eqref{eq:e-classes} in different fibres), so one
gets a map
\[
 \beta: BAut^0(M) \longrightarrow Conf_5^{ord}(\CP{1})/PSL_2(\C)
\]
unique up to homotopy, where $Aut^0(M)$ is the subgroup of
symplectic automorphisms acting trivially on homology. The
monodromy of the universal family of blowups over
$\BB^{ord}/PSL_3(\C)$ gives a map
\[
 \alpha: \BB^{ord}/PSL_3(\C) \longrightarrow BAut^0(M)\,
\]
such that $\beta \circ \alpha$ is homotopy equivalent to the
previous isomorphism of spaces. Hence $\alpha$ induces an
injective homomorphism from $\Gamma_5^{ord} =
\pi_1(\BB^{ord}/PSL_3(\C))$ to $\pi_0(Aut(M))$. By taking up the
discussion from the previous example, one can see that the image
of that homomorphism is generated by squared Dehn twists, hence
maps trivially to $\pi_0(\Diff(M))$. To take into account maps
which act nontrivially on homology, one should introduce an
extension $\Gamma_5$ of $\Gamma_5^{ord}$ by the Weyl group
$W(D_5)$, which is the automorphism group of the lattice $H_2(M)$
preserving $c_1$ and the intersection product (this is slightly
larger than the extension $\pi_1^{orb}(\BB/PSL_3(\C) \times S_5)$
that one gets from passing to unordered configurations). There is
a corresponding extended map $\Gamma_5 \rightarrow \pi_0(Aut(M))$,
which is obviously also injective, and which fits into a
commutative diagram
\begin{equation} \label{eq:d5-diagram}
\xymatrix{
 {\Gamma_5} \ar[d] \ar@{^{(}->}[rr] && {\pi_0(Aut(M))} \ar[d] \\
 {W(D_5)} \ar[rr] && {\pi_0(\Diff(M))}.
}
\end{equation}
\end{example}

$M = \CP{2} \# 5 \overline{\CP{2}}$ also occurs as space of
parabolic rank two odd degree bundles with fixed determinant and
weights $1/2$ on the five-pointed sphere (this gives another
explanation for the isomorphism of configuration spaces in Example
\ref{th:dp5}). In terms of flat connections, one can write it as
\begin{equation} \label{eq:flat}
 M = \{ A_1 \dots A_5 \in C_{1/2} \suchthat A_1 \dots A_5 = I \}/ PU(2).
\end{equation}
where $S^2 \iso C_{1/2} \subset SU(2)$ is the conjugacy class of
$diag(i,-i)$, and $PU(2)$ acts by simultaneous conjugation. There
is an obvious action of the mapping class group of the
five-pointed sphere on $M$. The gauge-theoretic definition shows
that the action is by symplectomorphisms, and up to symplectic
isotopy one can identify it with the top $\rightarrow$ from
\eqref{eq:d5-diagram} restricted to $\pi_1^{orb}(\BB/PSL_3(\C)
\times S_5) \subset \Gamma_5$. The injectivity of this map is
interesting because of the (conjectural) relation between
symplectic Floer homology and certain gauge theoretic invariants
of knots \cite{collin-steer99}. As a final remark, note that there
is a striking similarity between \eqref{eq:flat} and the
definition of the regular pentagon space: indeed, if one replaces
$C_{1/2}$ with the conjugation class of $diag(e^{\pi i
\alpha},e^{-\pi i \alpha})$ for some $\alpha < 2/5$, GIT arguments
show that the resulting space is symplectically deformation
equivalent to the pentagon space (as one passes the critical
weight $2/5$, the space undergoes a single blowup). This makes the
difference between the behaviour of squared Dehn twists even more
remarkable.

If one goes further to $k = 6$, where the blowup is a cubic
surface in $\CP{3}$, the situation becomes considerably more
complicated, mainly because the notion of general position
involves an additional condition on conics. Take the space
$\BB^{ord}/PSL_3(\C)$ of ordered configurations of points in
general position, which is the same as the moduli space of marked
cubic surfaces. A theorem of Allcock \cite{allcock02} says that
this is a $K(\Gamma_6^{ord},1)$, and the group $\Gamma_6^{ord}$ is
quite large: it contains infinitely generated normal subgroups
\cite{allcock-carlson-toledo99}. For purposes of comparison with
$\pi_0(Aut)$, the right group $\Gamma_6$ is an extension of
$\Gamma_6^{orb}$ by $W(E_6)$, which is the orbifold fundamental
group of the moduli space of cubic surfaces. Libgober
\cite{libgober78} proved that $\Gamma_6$ is a quotient of the
generalized braid group $B(E_6)$, and Looijenga \cite{looijenga98}
has given an explicit presentation of it. The last-mentioned paper
also contains a discussion of the $k = 7$ case, in which
$\Gamma_7$ is the orbifold fundamental group of the moduli space
of {\em non-hyperelliptic} genus three curves. We will return to
these del Pezzo surfaces in Example \ref{th:dp2} below.

\begin{example} \label{th:double-blowup}
Here is another, even simpler, application of Lemma
\ref{th:no-monodromy}. For any algebraic surface $\MM$, there is a
smooth family $\EE \rightarrow \BB$ over the configuration space
$\BB = Conf_2(\MM)$, whose fibre at $b$ is the blowup $Bl_b(\MM)$.
Take an ample line bundle $\Lambda$ on $\MM$, and equip each
blowup with $\LL_b = \Lambda^{\otimes d} \otimes
\OO(-E_1-E_2)^{\otimes e}$ for some $d \gg e \gg 0$. Both the
family and the line bundle extend to the compactification
$\overline{\BB} = Hilb_2(\MM)$ where the two points are allowed to
come together, and the fibres over the discriminant $\Delta$ have
ordinary double points. The monodromy around the meridian is a
Dehn twist along a Lagrangian sphere in the class $E_1 - E_2$, and
using the short exact sequence
\[
 1 \rightarrow \pi_1(\MM)^2 \rightarrow \pi_1(\BB)
 \rightarrow \Z/2 \rightarrow 1
\]
one sees that the square of this Dehn twist is isotopic to the
identity. This is actually a local phenomenon: $\C^2 \#
2\overline{\CP{2}}$, with the two exceptional divisors having equal
area, contains a Lagrangian sphere whose squared Dehn twist is
isotopic to the identity in the compactly supported symplectic
automorphism group. This can be implanted into $M \#
2\overline{\CP{2}}$ for any closed symplectic four-manifold $M$, as
long as the area of the exceptional divisors remains equal and
sufficiently small.
\end{example}

\section{Floer and quantum homology}

\subsection{}
Fix a closed symplectic four-manifold $M$ with $H^1(M;\R) = 0$, and a
coefficient field $\K$ ($\K = \Q$ will do in all the basic examples,
but including positive characteristic fields gives slightly sharper
general results). The universal Novikov field $\Lambda$ over $\K$ is
the field of formal series
\[
 f(q) = \sum_{d \in \R} a_d q^d
\]
with coefficients $a_d \in \K$, with the following one-sided growth
condition: for any $D \in \R$ there are at most finitely many $d \leq
D$ such that $a_d \neq 0$. Floer homology associates to any $\phi \in
Aut(M)$ a finite-dimensional $\Z/2$-graded $\Lambda$-vector space,
the Floer homology group
\[
HF_*(\phi) = HF_0(\phi) \oplus HF_1(\phi),
\]
and these groups come with the following additional structure:
\begin{itemize} \itemsep1em
\item
There is a distinguished element $e \in HF_0(id)$ and a distinguished
linear map $p: HF_0(id) \rightarrow \Lambda$.
\item
For any $\phi,\psi$ there is a canonical product, the so-called
pair-of-pants product
\[
\ast = \ast_{\phi,\psi}: HF_*(\phi) \otimes_\Lambda HF_*(\psi)
\longrightarrow HF_*(\phi\psi)
\]
\item
For any $\phi,\psi$ there is a conjugation isomorphism
\[
c_{\phi,\psi}: HF_*(\phi) \stackrel{\iso}{\longrightarrow} HF_*(\psi
\phi \psi^{-1})
\]
\item
For any smooth path $\lambda: [0;1] \rightarrow Aut(M)$ there is a
canonical continuation element $I_{\lambda} \in
HF_0(\lambda_0^{-1}\lambda_1)$.
\end{itemize}
We now write down a rather long list of axioms satisfied by Floer
homology theory. The aim is partly pedagogical, since this compares
unfavourably with the later formulation in terms of a topological
quantum field theory.
\begin{itemize} \itemsep1em
\item
$\ast$ is associative, in the sense that the two possible ways of
bracketing give the same trilinear map $HF_*(\phi) \otimes HF_*(\psi)
\otimes HF_*(\eta) \rightarrow HF_*(\phi\psi\eta)$. It is
commutative, which means that the following diagram commutes:
\[
\xymatrix{
 HF_*(\phi) \otimes HF_*(\psi) \ar[d]^{\ast}
 \ar[rrr]^-{\text{(signed) exchange}} &&&
 HF_*(\psi) \otimes HF_*(\phi) \ar[d]^{\ast} \\
 HF_*(\phi \psi) \ar[rrr]^-{c_{\phi\psi,\phi^{-1}}} &&&
 HF_*(\psi \phi)
}
\]
$e \in HF_*(id)$ is a two-sided unit for $\ast$, and for any $\phi$
we get a nondegenerate pairing between $HF_*(\phi)$ and
$HF_*(\phi^{-1})$ by setting $\leftsc x,y \rightsc = p(x \ast y)$.

\item
$c_{\phi,id}$ is the identity for any $\phi$, and so is
self-conjugation $c_{\phi,\phi}$ for any $\phi$. Conjugation
isomorphisms are well-behaved under composition,
$c_{\psi\phi\psi^{-1},\eta} \circ c_{\phi,\psi} = c_{\phi,\eta\psi}$.
They are compatible with pair-of-pants products, $c_{\phi,\eta}(x)
\ast c_{\psi,\eta}(y) = c_{\phi\psi,\eta}(x \ast y)$. Moreover,
conjugation $c_{id,\phi}: HF_*(id) \rightarrow HF_*(id)$ for any
$\phi$ leaves $e$ and $p$ invariant.

\item
Any constant path $\lambda$ gives rise to the element $I_\lambda = e
\in HF_*(id)$. Two paths which are homotopic rel endpoints have the
same continuation elements. Concatenation of paths corresponds to
product of continuation elements, $I_{\lambda \circ \mu} = I_\lambda
\ast I_\mu$. Next, if we compose a path $\lambda$ with a fixed map
$\phi$, more precisely if $(L_\phi\lambda)_t = \phi\lambda_t$ and
$(R_\phi\lambda)_t = \lambda_t\phi$, then
\[
 I_{L_\phi\lambda} = I_{\lambda}, \qquad
 I_{R_\phi\lambda} = c_{\lambda_0^{-1}\lambda_1,\phi^{-1}}(I_\lambda).
\]
\end{itemize}

\begin{remark}
We cannot pass this monument to abstract nonsense without lifting our
hat to gerbes. For simplicity we consider only finite cyclic gerbes,
so suppose that $X$ is a connected topological space carrying a
bundle of projective spaces $\CP{n} \rightarrow E \rightarrow X$ with
a $PU(n+1)$-connection, and $\Omega X$ the based loop space. To any
$\phi \in \Omega X$ one can associate the monodromy $m_{\phi} \in
PU(n+1)$. Take the set of all preimages of $m_\phi$ in $U(n+1)$, and
let $I(\phi)$ be the $\C$-vector space freely generated by this set.
(1) $I(\text{\it constant path})$ is the group ring $\C[\Z/(n+1)]$,
and we can define a canonical element $e$ and linear map $p$ as
usual. (2) Since $m_{\phi\psi} = m_{\phi}m_{\psi}$, multiplication in
$U(n+1)$ defines a composition map $I(\phi) \otimes I(\psi)
\rightarrow I(\phi\psi)$. (3) Conjugation with $m_\psi$ gives rise to
an isomorphism $I(\phi) \rightarrow I(\psi\phi\psi^{-1})$. (4) For
any homotopy $\lambda_t$ in $\Omega X$ one can define a preferred
element of $I(\lambda_1^{-1}\lambda_0)$ by deforming
$\lambda_1^{-1}\lambda_0$ to the constant path, and taking $e$ there.
This satisfies all the properties stated above.
\end{remark}

The first consequence of the axioms is that $HF_*(id)$ is a graded
commutative algebra with unit $e$. Actually, the trace $p$ makes it
into a Frobenius algebra. The conjugation maps $c_{id,\phi}$ define
an action of $Aut(M)$ on $HF_*(id)$ by Frobenius algebra
automorphisms, and this descends to an action of $\pi_0(Aut(M))$. To
see that, note that for any $x \in HF_*(id)$ and any path $\lambda$
starting at $\lambda_0 = id$, with corresponding reversed path
$\bar\lambda$, we have $c_{\lambda_1,\lambda_1} = id$ and
$I_{\lambda} \ast I_{\bar\lambda} = e$, hence
\[
 c_{id,\lambda_1}(x)
 = c_{id,\lambda_1}(x) \ast c_{\lambda_1,\lambda_1}(I_{\lambda})
 \ast I_{\bar\lambda} \\
 = c_{\lambda_1,\lambda_1}(x \ast I_{\lambda})
 \ast I_{\bar\lambda} \\
 = x \ast I_\lambda \ast I_{\bar\lambda} = x.
\]
$HF_*(id)$ acts on each $HF_*(\phi)$ by left pair-of-pants product
(one could equally use the product on the right, since for $x \in
HF_*(id)$ and $y \in HF_*(\phi)$, $y \ast x = (-1)^{deg(x)deg(y)}
c_{\phi,id}(x \ast y) = (-1)^{deg(x)deg(y)} x \ast y$). Here are some
simple properties of the module structure, directly derived from the
axioms:

\begin{lemma} \label{th:conjugation-action}
(1) $x,c_{id,\phi}(x) \in HF_*(id)$ act in the same way on $y \in
HF_*(\phi)$. (2) Up to isomorphism of $HF_*(id)$-modules,
$HF_*(\phi)$ is an invariant of $[\phi] \in \pi_0(Aut(M))$. (3) There
is a nondegenerate pairing $HF_*(\phi) \otimes HF_*(\phi^{-1})
\rightarrow \Lambda$ satisfying $\leftsc x \ast y, z \rightsc =
(-1)^{deg(x)deg(y)} \leftsc y, x \ast z \rightsc$ for all $x \in
HF_*(id)$, $y \in HF_*(\phi)$, $z \in HF_*(\phi^{-1})$.
\end{lemma}

\proof (1) $x \ast y = c_{\phi,\phi}(x \ast y) = c_{id,\phi}(x) \ast
c_{\phi,\phi}(y) = c_{id,\phi}(x) \ast y$. (2) For any path
$\lambda$, right multiplication with $I_\lambda$ is an isomorphism
$HF_*(\lambda_0) \rightarrow HF_*(\lambda_1)$ which commutes with
left multiplication by elements of $HF_*(id)$. (3) The pairing is
defined as $\leftsc y,z \rightsc = p(y \ast z)$, and obviously has
the desired properties. \qed

\subsection{}
By a theorem of Piunikhin-Salamon-Schwarz
\cite{piunikhin-salamon-schwarz94}, Ruan-Tian \cite{ruan-tian95}, and
Liu-Tian \cite{liu-tian01}, $HF_*(id)$ is canonically isomorphic to
the (small) {\em quantum homology ring} $QH_*(M)$. As a vector space,
this is simply $H_*(M;\Lambda)$ with the grading reduced to $\Z/2$.
The identity $e \in HF_*(id)$ is the fundamental class $[M]$, and the
linear map $p$ is induced from collapse $M \rightarrow point$. The
action of $\pi_0(Aut(M))$ is the obvious action of symplectomorphisms
on the homology of our manifold. The only non-topological element is
the quantum intersection product, which corresponds to the
pair-of-pants product in Floer cohomology, hence will be denoted by
the same symbol $\ast$. It is defined by
\[
 (x_0q^0 \ast y_0q^0) \cdot z_0q^0 =
 \sum_{A \in H_2(M;\Z)} \Phi_{3,A}(x_0,y_0,z_0) \, q^{\o(A)}
\]
for $x_0,y_0,z_0 \in H_*(M;\K)$, where $\cdot$ is the ordinary
intersection pairing with $\Lambda$-coefficients, and
$\Phi_{3,A}(x_0,y_0,z_0) \in \K$ the simplest kind of genus zero
Gromov invariant, counting pseudo-holomorphic spheres in class $A$
with three marked points lying on suitable representatives of
$x_0,y_0,z_0$ respectively. Note that since symplectic four-manifolds
are {\em weakly monotone}, we can (and will) use the older approach
of Ruan-Tian \cite{ruan-tian94} and McDuff-Salamon
\cite{mcduff-salamon} to define Gromov invariants with coefficients
in an arbitrary field $\K$. The leading term
$\Phi_{3,0}(x_0,y_0,z_0)q^0$ counting constant pseudo-holomorphic
curves is the ordinary triple intersection pairing, so the leading
term in the quantum product is the ordinary intersection product.

\begin{prop} \label{th:liu} {\rm (\cite[Corollary 1.6]{mcduff-salamon96b},
largely based on results of \cite{liu96})} Let $M$ be a closed
symplectic four-manifold which is minimal, and not rational or ruled.
Then $\Phi_{3,A} = 0$ for all $A \neq 0$. \qed
\end{prop}

Among the cases not covered by the Proposition, rational surfaces are
of primary interest because of the connection to classical
enumerative problems in projective geometry. Here is a very simple
example:

\begin{example} \label{th:qh-delpezzo}
We will be using some representation theory of finite groups, so let
$char(\K) = 0$ throughout the following computation. Take $M = \CP{2}
\# k \overline{\CP{2}}$, $5 \leq k \leq 8$, equipped with its
monotone symplectic structure, normalized to $[\o] = c_1$.
 Monotonicity simplifies the structure
of the quantum product considerably: in the expansion
\[
 x \ast y = (x \cap y) + (x \ast_1 y)q + (x \ast_2 y)q^2 + \dots
\]
the $q^d$ term has degree $2d-4$ with respect to the ordinary grading
of $H_*(M;\Lambda)$, in particular the terms $q^5,q^6,\dots$ all
disappear. Fixing some compatible (integrable) complex structure, one
finds that the only holomorphic spheres with $c_1(A) = 1$ are the
exceptional divisors, of which there is precisely one for each
element of $\EE = \{A \in H_2(M;\Z) \suchthat c_1(A) = 1, \; A \cdot
A = -1\}$. By the divisor axiom for Gromov invariants, these classes
$A$ satisfy $\Phi_{3,A}(x,y,z) = (x \cdot A) (y \cdot A) (z \cdot A)$
for $x,y,z \in H_2(M;\K)$, and hence
\begin{equation} \label{eq:ast1}
 x \ast_1 y = \sum_{A \in \EE} (x \cdot A) (y \cdot A) A.
\end{equation}

Let $K$ be the Poincar{\'e} dual of $-c_1$, and $K^\perp \subset
H_2(M;\K)$ its orthogonal complement with respect to the intersection
form. {\bf Fact:} {\em for all $x,y \in K^\perp$, $x \ast_1 y$ is a
multiple of $K$}. This follows from \eqref{eq:ast1} by explicit
computation \cite[Proposition 3.5.5]{bayer-manin01}. For the most
complicated cases $k = 7,8$ one can also use a trick from \cite[p.\
33]{demazure-pinkham-teissier77}: $\bar{A} = (k-6)(A \cdot K) K - A$
is an involution of $H_2(M;\Z)$ preserving $K$ and the intersection
form. It acts freely on $\EE$, and the contributions of $A$ and
$\bar{A}$ to $x \ast_1 y$ add up to a multiple of $K$. Next, let $W$
be the group of linear automorphisms of $H_2(M;\Z)$ which preserve
the intersection form, and leave $K$ fixed. This is a reflection
group of type $D_5,E_6,E_7$ or $E_8$ and it acts irreducibly on
$K^\perp$. Moreover, each element of $W$ can be realized by a
symplectic automorphism of $M$, and so the quantum product is
$W$-equivariant (this can also be checked by a direct computation of
Gromov invariants, without appealing to the $Aut(M)$-action).
Therefore, both $\ast_1: (K^\perp)^{\otimes 2} \rightarrow \K K
\subset H_2(M;\K)$ and $\ast_2: (K^\perp)^{\otimes 2} \rightarrow
H_4(M;\K) = \K$ must be scalar multiples of the intersection form. We
record this for later use, {\bf Fact:} {\em There is a $z \in
QH_*(M)$ of the form $z = [point] + \alpha_1 K q + \alpha_2 [M]q^2$
for some $\alpha_1,\alpha_2 \in \K$, such that for all $x,y \in
K^\perp$, $x \ast y = (x \cdot y) z$.}
\end{example}

\subsection{}
The case $\phi = id$ is misleading in so far as for a general
symplectic automorphism $\phi$, $HF_*(\phi)$ has no known
interpretation in terms of topology or Gromov-Witten invariants, and
is hard or impossible to compute. Our insight into Dehn twists and
their squares depends entirely on the following result:

\begin{prop} \label{th:exact}
For any Lagrangian sphere $L \subset M$, there is a long exact
sequence
\[
 \xymatrix{
 H_*(S^2;\Lambda) \ar[r] &
 QH_*(M) \ar[r]^-{G} &
 HF_*(\tau_L) \ar@/^1pc/[ll]^-{\partial}
 }
\]
where the grading of $H_*(S^2;\Lambda)$ is reduced to a
$\Z/2$-grading, $\partial$ has odd degree, and $G$ is a map of
$QH_*(M)$-modules.
\end{prop}

The origins of this will be discussed extensively later, but for now
let's pass directly to applications. Let $I_l \subset QH_*(M)$ be the
ideal generated by $l = [L]q^0$.

\begin{lemma} \label{th:ideal}
$dim_\Lambda I_l = 2$, and moreover $I_l$ is contained in $QH_0(M)$.
\end{lemma}

\proof Assume first that $char(\K) \neq 2$. Since $L \cdot L =
-2$, we know that $l$ is nontrivial and linearly independently
from $l \ast l = -2[point] + \dots$, so $dim_\Lambda I_l \geq 2$.
The other half uses the Picard-Lefschetz formula \eqref{eq:p-l}.
Since $(\tau_L)_*(l) = -l$, multiplication with $l$ is an
endomorphism of $QH_*(M)$ which exchanges the $\pm 1$ eigenspaces
of $(\tau_L)_*$. The $+1$ eigenspace has codimension one, and the
$-1$ eigenspace has dimension one, and so the kernel of the
multiplication map has codimension at most two, which means that
its image has dimension at most two.

Without assumptions on the characteristic, one has to argue slightly
more carefully as follows. We know that $[L] \in H_2(M;\Z)$ is
nontrivial and primitive, so there is a $w \in H_2(M;\Z)$ with $w
\cdot [L] = 1$. Denote the induced element of $H_2(M;\K)$ equally by
$w$. Then $l \ast w = [point] + \dots$, from which it follows as
before that $dim_\Lambda I_l \geq 2$. From the Picard-Lefschetz
formula one gets
\[
   w \ast l + ((w \ast l) \cdot l) l
 = (\tau_L)_*(w \ast l)
 = (\tau_L)_*(w) \ast (\tau_L)_*(l)
 = -w \ast l - l \ast l,
\]
which shows that $l \ast l$ lies in the linear subspace generated by
$l$ and $w \ast l$; and similarly for any $x \in QH_*(M)$,
\[
 w \ast x + ((w \ast x) \cdot l) l =
 (\tau_L)_*(w \ast x) = w \ast x + l \ast x + (x \cdot l)(w \ast l + l \ast l)
\]
which shows that $l \ast x$ lies in the subspace generated by $l$ and
$w \ast l$. \qed

\begin{lemma}
The kernel of any $QH_*(M)$-module map $G: QH_*(M) \rightarrow
HF_*(\tau_L)$ must contain $I_l$.
\end{lemma}

\proof Let $w$ be as in the proof of the previous Lemma. From Lemma
\ref{th:conjugation-action}(1) we know that for any $y \in
HF_*(\tau_L)$, $l \ast y = (\tau_L)_*(w) \ast y - w \ast y = 0$.
Hence $G(l) = G(l \ast e) = l \ast G(e) = 0$, and therefore also $G(x
\ast l) = 0$ for any $x$. \qed

For the long exact sequence from Proposition \ref{th:exact}, this
means that the kernel of $G$ is precisely $I_l$ and that the
differential $\delta$ is zero, showing that
\[
 HF_*(\tau_L) \iso QH_*(M)/I_l
\]
as a $QH_*(M)$-module. Now suppose that $\tau_L^2$ is symplectically
isotopic to the identity. By Lemma \ref{th:conjugation-action}(2) we
have an isomorphism $HF_*(\tau_L^{-1}) \iso HF_*(\tau_L)$ of
$QH_*(M)$-modules, and part (3) of the same Lemma shows that there is
a nondegenerate pairing on $QH_*(M)/I_l$ which satisfies $\leftsc x
\ast y, z \rightsc = \pm\leftsc y, x \ast z \rightsc$. Taking $y = e$
shows that $\leftsc x, z \rightsc = \leftsc e, x \ast z \rightsc$, so
the pairing comes from the linear map $\leftsc e, - \rightsc$ and the
quantum product on $QH_*(M)/I_l$.

\begin{cor} \label{th:w}
If $\tau_L^2$ is symplectically isotopic to the identity, the
quotient algebra $QH_*(M)/I_l$ is Frobenius. In particular, any
linear subspace $W \subset QH_0(M)/I_l$ which satisfies $x \cdot y =
0$ for all $x,y \in W$ must satisfy $dim_\Lambda W \leq \half
dim_\Lambda QH_0(M)/I_l$.
\end{cor}

The first part is just the outcome of the preceding discussion, and
the second part is an elementary fact about Frobenius algebras: $W$
is an isotropic subspace with respect to the pairing, whence the
bound on the dimension.

\begin{cor} \label{th:minimal}
Let $M$ be a closed minimal symplectic four-manifold with
$H^1(M;\R) = 0$, and not rational or ruled. Suppose that
$dim\,H_2(M;\K) \geq 3$. Then for every Lagrangian sphere $L
\subset M$, $\tau_L^2$ is not symplectically isotopic to the
identity, hence fragile.
\end{cor}

\proof From Proposition \ref{th:liu},
\[
QH_*(M)/I_l = H_*(M;\Lambda)/(\Lambda l \oplus \Lambda [point])
\]
with the algebra structure induced by the ordinary intersection
product. In particular, $W = H_2(M;\Lambda)/\Lambda l$ is a
subspace satisfying the conditions of Corollary \ref{th:w}, and
$dim\, W = dim\, H_2(M;\Lambda) - 1 > \half dim\, H_2(M;\Lambda) =
\half dim\, QH_*(M)/I_l$. \qed

Example \ref{th:double-blowup} shows that the minimality assumption
cannot be removed. The condition that $M$ should not be rational
excludes the case of $S^2 \times S^2$ discussed in Example
\ref{th:s2s2}. As for the final assumption $dim\, H_2(M;\K) \geq 3$,
a lack of suitable examples makes it hard to decide whether it is
strictly necessary. In the algebro-geometric world, there are minimal
surfaces of general type with Betti numbers $b_1(M) = 0$, $b_2(M) =
2$ exist, but the Miyaoka inequality $\chi - 3 \sigma \geq
\frac{9}{2}\#\{nodes\}$ \cite{miyaoka84} implies that they do not
admit degenerations to nodal ones, thereby barring the main route to
constructing Lagrangian spheres in them. Moreover, the most common
explicit examples in the literature are uniformized by a polydisc, so
they cannot contain any embedded spheres with nonzero
selfintersection.

\begin{example} \label{th:dp2}
Take $M = \CP{2} \# k\overline{\CP{2}}$, $5 \leq k \leq 8$, with a
monotone symplectic form. As in Example \ref{th:qh-delpezzo} we use a
coefficient field with $char(\K) = 0$. The computation carried out
there shows that for any $x,y \in K^\perp$, $x \ast y = -\half(x
\cdot y) l \ast l \in I_l$. Hence, the image of $K^\perp$ in
$QH_*(M)/I_l$, which is of dimension
\begin{equation} \label{eq:k2}
k - 1 > \half(k+1) = \half dim\, QH_*(M)/I_l,
\end{equation}
violates the conditions of Corollary \ref{th:w}. It follows that
in contast with the situation for $k \leq 4$, squared Dehn twists
are never symplectically isotopic to the identity. For $k = 5$, we
already saw some cases of this phenomenon in Example \ref{th:dp5},
and as explained there, this goes well with the intuition provided
by the topology of moduli spaces. In a slightly different
direction, one should note that the nontriviality of $\tau^2$ has
implications for the $\pi_1$ of spaces of symplectic embeddings of
$k$ balls into $\CP{2}$, via the symplectic interpretation of
blowup, see e.g.\ \cite{biran96}.
\end{example}

It would be interesting to extend the entire discussion to arbitrary
(not monotone) symplectic forms on rational four-manifolds. Although
the Gromov invariants are constant under deformations of the
symplectic class, the exponents $q^{\o(A)}$ change, which affects the
algebraic structure of the quantum homology ring, and thereby the
criterion which we have used to explore the nature of squared Dehn
twists. As a sample question, take a Lagrangian sphere $L$ on, say,
the cubic surface, and then perturb the symplectic class in a generic
way subject only to the condition that $L$ continues to be
Lagrangian. Is it true that then, $QH_*(M)/I_l$ becomes semisimple?
This is relevant because semisimple algebras are obviously Frobenius
(see \cite{bayer-manin01} for a proof of the generic semisimplicity
of $QH_*(M)$ itself).

Finally, we turn to the proof of Theorem \ref{th:one} stated in the
introduction (together with Corollary \ref{th:fragility}, this also
proves Corollary \ref{th:two}). Let $M \subset \CP{n+2}$ be a
nontrivial complete intersection of degrees ${\bf d} =
(d_1,\dots,d_n)$, $n \geq 1$ and $d_k \geq 2$, with the symplectic
structure $\o$ induced by the Fubini-Study form $\o_{FS}$, which we
normalize to $\o_{FS}^{n+2} = 1$. Each such $M$ contains a Lagrangian
sphere, which can be obtained as vanishing cycle in a generic pencil
of complete intersections. Moreover,
\begin{align*}
 \pi_1(M) & = 1, \\
 \chi(M) & =
 \frac{1}{2} \Big(\prod_k d_k\Big) \Big[ \Big( \sum_k d_k - (n+3)
 \Big)^{\!2} + \sum_k d_k^2 - (n+3) \Big], \\
  c_1(M) & = \Big(n+3-\sum_k d_k\Big)[\o].
\end{align*}
With the exception of six choices of degrees ${\bf d} = (2), (3),
(4), (2,2), (2,3), (2,2,2)$, $c_1(M)$ is a negative multiple of
$[\o]$, so $M$ is minimal and of general type, and $\chi(M) > \sum_k
d_k(d_k-1) \geq 6$, which means $b_2(M) \geq 4$, so Corollary
\ref{th:minimal} applies. Out of the remaning cases, three are $K3$
surfaces, ${\bf d} = (4), (2,3), (2,2,2)$, to which Corollary
\ref{th:minimal} also applies. The other three are ${\bf d} = (2)$
which is the quadric $\CP{1} \times \CP{1}$, hence excluded from the
statement of Theorem \ref{th:one}, and ${\bf d} = (3), (2,2)$ which
are the del Pezzo surfaces of rank $k = 6,5$ respectively, and
therefore fall under Example \ref{th:dp2}.

\subsection{}
As promised, we will now present Floer homology theory as a TQFT in
$1+1$ dimensions ``coupled with'' symplectic fibrations. This is a
generalization of the setup from \cite{piunikhin-salamon-schwarz94}
where only the trivial fibration was allowed (Lalonde has recently
introduced a very similar generalization, but his intended
applications are quite different). Throughout the following
discussion, all symplectic fibrations have fibres isomorphic to $M$,
without any specific choice of isomorphism. The basic data are
\begin{itemize}
\item \itemsep1em
For any symplectic fibration $F \rightarrow Z$ over an oriented
circle $Z$, we have a Floer homology group $HF_*(Z,F)$.
\item
For any isomorphism $\Gamma: F_1 \rightarrow F_2$ between such
fibrations covering an orientation-preserving diffeomorphism $\gamma:
Z_1 \rightarrow Z_2$, there is an induced canonical isomorphism
$C(\gamma,\Gamma): HF_*(Z_1,F_1) \rightarrow HF_*(Z_2,F_2)$.
\item
Let $S$ be a connected compact oriented surface with $p+q$ boundary
circles. We arbitrarily divide the circles into positive and negative
ones, and reverse the natural induced orientation of the latter, so
that
\[
\partial S = \bar{Z}^-_1 \cup \dots \cup \bar{Z}^-_p \cup Z^+_{p+1}
\cup \dots \cup Z^+_{p+q}.
\]
Given a symplectic fibration $E \rightarrow S$, with restrictions
$F^\pm_k = E|Z^\pm_k$, we have a relative Gromov invariant
\[
G(S,E): \bigotimes_{k=1}^p HF_*(Z^-_k,F^-_k) \longrightarrow
\bigotimes_{k=p+1}^{p+q} HF_*(Z^+_k,F^+_k).
\]
This is independent of the way in which the $Z^-,Z^+$ are numbered,
up to the usual signed interchange of factors in the tensor product.
\end{itemize}
The maps $C(\gamma,\Gamma)$, sometimes omitted from more summary
expositions, are a natural part of the theory: after all, the
``cobordism category'' is more properly a 2-category
\cite{tillmann98}, and the algebraic framework should reflect this.
The TQFT axioms are
\begin{itemize} \itemsep1em
\item
The identity automorphism of each $(Z,F)$ induces the identity
$C(id_Z,id_F)$ on Floer homology. The maps $C(\gamma,\Gamma)$ are
well-behaved under composition of isomorphisms. Moreover, if
$(\gamma^t,\Gamma^t)$, $t \in [0;1]$ is a smooth family of
isomorphisms $F_1 \rightarrow F_2$, then $C(\gamma^0,\Gamma^0) =
C(\gamma^1,\Gamma^1)$.

\item
Let $\xi: S_1 \rightarrow S_2$ be an orientation-preserving
diffeomorphism of surfaces, which respects the decomposition of the
boundary into positive and negative circles, and suppose that this is
covered by an isomorphism $\Xi: E_1 \rightarrow E_2$ of symplectic
fibrations. Let $\gamma^\pm_k$, $\Gamma^{\pm}_k$ be the restriction
of $\xi,\Xi$ to the boundary components. Then the following diagram
commutes:
\[
\xymatrix{
 {\bigotimes_{k=1}^p HF_*(Z_{1,k}^-,F_{1,k}^-)}
 \ar[rr]^-{G(S_1,E_1)}
 \ar[d]_-{\bigotimes_{k=1}^p C(\gamma_k^-,\Gamma_k^-)}
 &&
 {\bigotimes_{k=p+1}^{p+q} HF_*(Z_{1,k}^+,F_{1,k}^+)}
 \ar[d]^-{\bigotimes_{k=p+1}^{p+q} C(\gamma_k^+,\Gamma_k^+)}
 \\
 {\bigotimes_{k=1}^p HF_*(Z_{2,k}^-,F_{2,k}^-)}
 \ar[rr]^-{G(S_2,E_2)}
 &&
 {\bigotimes_{k=p+1}^{p+q} HF_*(Z_{2,k}^+,F_{2,k}^+)}
}
\]

\item
Take any $F \rightarrow Z$ and pull it back by projection to a
fibration $E \rightarrow S = [1;2] \times Z$. (a) If we take $Z_1 =
\{1\} \times Z$ negative and $Z_2 = \{2\} \times Z$ positive, the
relative Gromov invariant $HF_*(Z,F) \rightarrow HF_*(Z,F)$ is the
identity. (b) Take both $Z_1,Z_2$ to be negative. Then the relative
Gromov invariant $HF_*(Z,F) \otimes HF_*(\bar{Z},\bar{F}) \rightarrow
\Lambda$, where $(\bar{Z},\bar{F})$ denotes orientation-reversal on
the base, is a nondegenerate pairing.

\item
The gluing or cut-and-paste axiom. Let $S_1,S_2$ be two surfaces
carrying symplectic fibrations $E_1,E_2$, and suppose that we have an
isomorphism $(\gamma,\Gamma)$ between the induced fibrations over the
$m$-th positive boundary circle of $S_1$ and the $n$-th negative one
of $S_2$. One can glue together the two boundary components to form a
surface $S = S_1 \cup S_2$ and a symplectic fibration $E$ over it,
and the associated relative invariant $G(S,E)$ is the composition
\[
\xymatrix{
 \Big(\bigotimes_k HF_*(Z_{1,k}^-,F_{1,k}^-)\Big) \otimes
 \Big(\bigotimes_{l \neq n} HF_*(Z_{2,l}^-,F_{2,l}^-)\Big)
 \ar[d]^-{G(S_1,E_2) \otimes id} \\
 \Big(\bigotimes_k HF_*(Z_{1,k}^+,F_{1,k}^+)\Big) \otimes
 \Big(\bigotimes_{l \neq n} HF_*(Z_{2,l}^-,F_{2,l}^-)\Big)
 \ar[d]_-{\iso}^-{\text{exchange and $C(\gamma,\Gamma)$}} \\
 \Big(\bigotimes_{k \neq m} HF_*(Z_{1,k}^+,F_{1,k}^+)\Big)
 \otimes \Big(\bigotimes_l HF_*(Z_{2,l}^-,F_{2,l}^-)\Big)
 \ar[d]^{id \otimes G(S_2,E_2)} \\
 \Big(\bigotimes_{k \neq m} HF_*(Z_{1,k}^+,F_{1,k}^+)\Big)
 \otimes \Big(\bigotimes_l HF_*(Z_{2,l}^+,F_{2,l}^+)\Big)
}
\]
\end{itemize}

How does this set of axioms for Floer homology imply the
previously used one? To any $\phi \in Aut(M)$ one can associate
the mapping torus $F_\phi = \R \times M / (t,x) \sim
(t-1,\phi(x))$, which is naturally a symplectic fibration over
$S^1 = \R/\Z$. Set $HF_*(\phi) = HF_*(S^1,F_\phi)$. Given an
isotopy $(\lambda_t)$ in $Aut(M)$ (with $\lambda_t$ constant for
$t$ close to the endpoints $0,1$), one can define an isomorphism
$\Gamma_\lambda: F_{\lambda_0} \rightarrow F_{\lambda_1}$ by
$[0;1] \times M \rightarrow [0;1] \times M$, $(t,x) \mapsto
(t,\lambda_t^{-1}\lambda_0(x))$, and the corresponding map
$C(id_{S^1},\Gamma_\lambda)$ is our previous $I_\lambda$. For any
$\psi$ there is a canonical isomorphism $\Gamma_{\phi,\psi}:
F_\phi \rightarrow F_{\psi\phi\psi^{-1}}$, $(t,x) \mapsto
(t,\psi(x))$, and we correspondingly define the conjugation
isomorphism $c_{\phi,\psi} = C(id_{S^1},\Gamma_{\phi,\psi})$. In
the case where $\phi = \psi$, the isomorphism $\Gamma_{\phi,\phi}$
can be deformed to the identity by rotating the base once, $(t,x)
\mapsto (t-\tau,\phi(x))$, and this explains the previously stated
property that $c_{\phi,\phi} = id$. Extracting the remaining
structure, such as the pair-of-pants and its properties, is staple
TQFT fare, which can be found in any expository account such as
\cite{segal00}.

Having come this far, we can make a straightforward extension to
the formalism, which is to replace symplectic fibrations by
Lefschetz fibrations in the sense of Definition
\ref{th:lefschetz}. This requires some small modifications of the
axioms, since even in the absence of critical points, our
definition of Lefschetz fibrations contains more information (the
two-form $\Omega$ on the total space) than that of symplectic
fibration. Essentially, one has to add another property saying
that relative Gromov invariants are unchanged under deformation of
a Lefschetz fibration. But we have spent enough time with
exercises in axiomatics, so we leave the precise formulation to
the reader. The essential new ingredient that comes from Lefschetz
fibrations is this: given any Lagrangian sphere $L \subset M$, one
can construct a Lefschetz fibration $E'$ over a disc $S'$ with a
single critical point, whose associated vanishing cycle is $L
\subset M$. By the Picard-Lefschetz theorem, the monodromy around
the boundary is isotopic to the Dehn twist $\tau_L \in Aut(M)$,
and the associated relative Gromov invariant provides a
distinguished element
\begin{equation} \label{eq:char}
\theta_L \stackrel{\text{def}}{=} G(S',E') \in HF_*(\tau_L).
\end{equation}
One can show that $E'$ is unique up to deformation, so this class is
independent of the details of the construction. Pair-of-pants product
with $\theta_L$ yields for any $\phi \in Aut(M)$ a canonical
homomorphism
\[
 HF_*(\phi) \longrightarrow HF_*(\phi \circ \tau_L),
\]
and the special case $\phi = id_M$ is the map $G$ from Proposition
\ref{th:exact}. The fact that this map is a homomorphism of
$QH_*(M)$-modules follows from associativity of the pair-of-pants
product.

\begin{remark} \label{th:general-exact}
One expects that Proposition \ref{th:exact} is a special case of a
more general long exact sequence, of the form
\begin{equation} \label{eq:general-exact}
 \xymatrix{
 HF_*(L,\phi(L)) \ar[r] &
 HF_*(\phi) \ar[r] &
 HF_*(\phi \circ \tau_L) \ar@/^1pc/[ll]
 }
\end{equation}
for an arbitrary $\phi \in Aut(M)$. The appearance of Lagrangian
intersection Floer homology means that in order to understand this
sequence, our framework should be further extended to an
``open-closed'' string theory, where the symplectic fibrations are
allowed to carry Lagrangian boundary conditions, see \cite{seidel01}.
In cases where Lagrangian Floer homology is well-behaved, such as
when $M$ is an exact symplectic manifold with boundary, the sequence
\eqref{eq:general-exact} can be readily proved by adapting of the
argument from \cite{seidel01}. When $M$ is a closed four-manifold,
$HF_*(L_0,L_1)$ is not always defined, but this should not be an
issue for the case for the group in \eqref{eq:general-exact}, since
the obstructions in the sense of \cite{oh95b, fooo} coming from $L$
and $\phi(L)$ ought to cancel out. With this in mind, the proof
should go through much as before, but there are still some technical
points to be cleared up, so we will stop short of claiming it as a
theorem.
\end{remark}

\begin{remark} \label{th:hamiltonian}
The condition $H^1(M;\R) = 0$ can be removed from the whole section,
at the cost of replacing $Aut(M)$ by $Aut^h(M)$. The Dehn twist along
a Lagrangian sphere is unique up to Hamiltonian isotopy; the axioms
for Floer homology remain the same except that the elements
$I_\lambda$ exist only for Hamiltonian isotopies; and since the basic
Proposition \ref{th:exact} continues to hold, so do all its
consequences. In the construction of the TQFT, one has to replace
symplectic fibrations by Hamiltonian fibrations (Lefschetz fibrations
as we defined them are already Hamiltonian).

A more interesting question is whether for $H^1(M;\R) \neq 0$, it can
happen that $\tau_L^2$ is symplectically and not Hamiltonian isotopic
to the identity. Assuming some unproved but quite likely statements,
one can give a negative answer to this at least in the case when $c_1
= \lambda[\o]$ for $\lambda<0$ and the divisibility of $c_1$ is $N
\geq 2$. Suppose that $\phi = \tau_L^{-2}$ is symplectically, but not
Hamiltonian, isotopic to the identity. There is a theorem of L\^e-Ono
\cite{le-ono95} which determines $HF_*(\phi)$ in the opposite sign
case where $\lambda>0$. It seems reasonable to expect this to hold in
our case too, so that
\begin{equation} \label{eq:le-ono}
 HF_*(\phi) \iso H_*(M;\underline{\Lambda})
\end{equation}
where $\underline{\Lambda}$ is a nontrivial local system of
$\Lambda$-coefficients determined by the flux of $\phi$. We also
assume the long exact sequence \eqref{eq:general-exact} for $\phi$,
which would be
\begin{equation} \label{eq:nonham-exact}
\xymatrix{
 H_*(L;\Lambda) \ar[r] &
 HF_*(\phi) \ar[r] &
 HF_*(\phi \circ \tau_L) \ar@/^1pc/[ll]^-{\delta \text{ (of degree $-1$)}}
}
\end{equation}
Floer homology groups are now $\Z/2N$-graded, and by combining
\eqref{eq:le-ono} with \eqref{eq:nonham-exact} and standard facts
about Novikov homology, one sees that $HF_*(\phi \circ \tau_L)$ is
concentrated in three adjacent degrees. On the other hand, we still
have $HF_*(\tau_L) \iso QH_*(M)/I_l$, which is nonzero in four
degrees, hence $HF_*(\tau_L^{-1}) \not\iso HF_*(\phi \circ \tau_L)$.
\end{remark}

\section{Pseudo-holomorphic sections and curvature}

\subsection{}
The aim of this section is to explain the proof of Proposition
\ref{th:exact}, but we start on a much more basic level with the
definition of Floer homology according to Hofer-Salamon
\cite{hofer-salamon95}, recast in fibre bundle language. Let
\[
p: F \longrightarrow S^1 = \R/\Z
\]
be a smooth proper fibration with four-dimensional fibres, equipped
with a closed two-form $\Omega$ whose restriction to each fibre is
symplectic. We have the corresponding symplectic connection $TF =
TF^h \oplus TF^v$. Let $\SS(S^1,F)$ be the space of all smooth
sections of $p$, and $\SS^h(S^1,F)$ the subspace of horizontal
sections, which are those with $d\sigma/dt \in TF^h$. To $\sigma \in
\SS^h(S^1,F)$ one can associate a linear connection $\nabla^\sigma$
on the pullback bundle $\sigma^*TF^v$,
\[
 \nabla^\sigma_{\partial_t} X = [d\sigma/dt,X].
\]
We say that $F$ has {\em nondegenerate horizontal sections} if for
every $\sigma$, there are no nonzero solutions of $\nabla^\sigma X =
0$. We will assume from now on that this is the case; then
$\SS^h(S^1,F)$ is finite, and one defines the Floer chain group as
\[
 CF_*(S^1,F) = \bigoplus_{\sigma \in \SS^h(S^1,F)} \Lambda \gen{\sigma}.
\]
The $\Z/2$ degree of a generator $\gen{\sigma}$ is determined by the
sign of $\det(1-R^\sigma)$, where $R^\sigma$ is the monodromy of
$\nabla^\sigma$ around $S^1$. There is a closed action one-form $da$
on $\SS(S^1,F)$ whose critical point set is $\SS^h(S^1,F)$, namely
\[
da(\sigma)X = \int_{S^1} \Omega(d\sigma/dt,X)\,dt
\]
for $X \in T_\sigma\SS(S^1,F) = \smooth(\sigma^*TF^v)$. Nondegeneracy
of the horizontal sections corresponds to the Morse nondegeneracy of
a local primitive $a$. Now take a smooth family of
$\Omega|F_t$-compatible almost complex structures $J_{F_t}$ on the
fibres. The negative gradient flow lines of $da$ with respect to the
resulting $L^2$ metric are the solutions of Floer's equation. In view
of later developments, we find it convenient to write the equation as
follows. Take $\pi = id_\R \times p: E = \R \times F \rightarrow S =
\R \times S^1$. Equip $S$ with its standard complex structure $j$,
and $E$ with the almost complex structure $J$ characterized by the
following properties: (1) $\pi$ is $(J,j)$-holomorphic; (2) The
restriction of $J$ to any fibre $E_{s,t}$ is equal to $J_{F_t}$; (3)
$J$ preserves the splitting of $TE$ into horizontal and vertical
parts induced by the pullback of $\Omega$. Then Floer's equation
translates into the pseudo-holomorphic section equation
\begin{equation} \label{eq:floer}
\begin{cases}
 \!\!\!\!\! & u: S \longrightarrow E, \quad \pi \circ u = id_S \\
 \!\!\!\!\! & Du \circ j = J \circ Du, \\
 \!\!\!\!\! & \lim_{s \rightarrow \pm\infty} u(s,\cdot) = \sigma_{\pm},
\end{cases}
\end{equation}
where $\sigma_{\pm} \in \SS^h(S^1,F)$. Note that for any $\sigma \in
\SS^h(S^1,F)$ there is a trivial or stationary solution $u(s,t) =
(s,\sigma(t))$ of \eqref{eq:floer}. We denote by
$\MM^*(S,E;\sigma_-,\sigma_+)$ the space of all {\em other}
solutions, divided by the free $\R$-action of translation in
$s$-direction; and by $\MM_0^*(S,E;\sigma_-,\sigma_+)$ the subspace
of those solutions whose virtual dimension is equal to zero. The
Floer differential on $CF_*(S^1,F)$ is defined by
\[
 \partial \gen{\sigma_-} \;\; = \!\!\!\! \sum_{
 \substack{
 \sigma_+ \in \SS^h(S^1,F) \\
 \!\!\!\! u \in \MM^*_0(S,E;\sigma_-,\sigma_+)\!\! } }
 \pm q^{\epsilon(u)} \gen{\sigma_+}
\]
where the energy is $\epsilon(u) = \int_S u^*\Omega \in (0;\infty)$,
and the sign is determined by coherent orientations which we will not
explain further. For this to actually work and give the correct Floer
homology $HF_* = H_*(CF_*,\partial)$, the $J_{F_t}$ need to satisfy a
number of generic ``transversality'' properties:
\begin{itemize}
\item
There are no non-constant $J_{F_t}$-holomorphic spheres of Chern
number $\leq 0$;
\item
If $v: S^2 \rightarrow F_t$ is a $J_{F_t}$-holomorphic map with Chern
number one, the image of $v$ is disjoint from $\sigma(t)$ for all
$\sigma \in \SS^h(S^1,F)$.
\item
The linearized operator $D\bar\partial_u$ attached to any solution of
Floer's equation is onto. This means that the spaces
$\MM^*(S,E;\sigma_-,\sigma_+)$ are all smooth of the expected
dimension.
\end{itemize}
The space of pseudo-holomorphic spheres with Chern number $\leq 0$
in a four-manifold has virtual dimension $\leq -2$, so that even
in a one-parameter family of manifolds the virtual dimension
remains negative. As for the images of pseudo-holomorphic spheres
with Chern number one, they form a codimension 2 subset in a
four-manifold, and the same thing holds in a family, so they
should typically avoid the image of any fixed finite set of
sections, which is one-dimensional. In both cases, the fact that
the condition is actually generic is proved by appealing to the
theory of somewhere injective pseudo-holomorphic curves, see for
instance \cite{mcduff-salamon}. The last requirement is slightly
more tricky because of the $\R$-symmetry on the moduli space; see
\cite{floer-hofer-salamon94} for a proof.

We now introduce the second ingredient of the TQFT, the relative
Gromov invariants. As a basic technical point, the surfaces with
boundary which we used to state the axioms must be replaced by
noncompact surfaces with a boundary at infinity. For ease of
formulation, we will consider only the case of $S = \R \times S^1$,
which is the one relevant for our applications. Let $\pi: E
\rightarrow S$ be a smooth proper fibration with four-dimensional
fibres, equipped with a closed two-form $\Omega$ whose restriction to
any fibre is symplectic. The behaviour of $E$ over the two ends of
our surface is governed by the following ``tubular ends''
assumptions: there are fibrations $p^\pm: F^\pm \rightarrow S^1$ with
two-forms $\Omega^{\pm}$ as before, with the property that the
horizontal sections are nondegenerate, and fibered diffeomorphisms
$\Psi^-: E|(-\infty;s_-] \longrightarrow (-\infty;s_-] \times F^-$,
$\Psi^+: E|[s_+;\infty) \longrightarrow [s_+;\infty) \times F^+$ for
some $s_-<s_+$, such that $(\Psi^{\pm})^*\Omega^{\pm} = \Omega$.

Take a positively oriented complex structure $j$ on $S$. We say that
an almost complex structure $J$ on $E$ is {\em semi-compatible} with
$\Omega$ if $\pi$ is $(J,j)$-holomorphic, and the restriction of $J$
to each fibre is compatible with the symplectic form in the usual
sense. With respect to the splitting $TE_x = TE_x^h \oplus TE_x^v$,
this means that
\begin{equation} \label{eq:jj}
J = \begin{pmatrix} j & 0 \\ J^{vh} & J^{vv} \end{pmatrix}
\end{equation}
where $J^{vv}$ is a family of compatible almost complex structures on
the fibres, and $J^{vh}$ is a $\C$-antilinear map $TE^h \rightarrow
TE^v$ (this corresponds to the ``inhomogeneous term'' in the theory
of pseudoholomorphic maps). We also need to impose some conditions at
infinity. Choose families of almost complex structures $J_{F^-_t}$,
$J_{F^+_t}$ on the fibres of $F^\pm$ which are admissible for Floer
theory, meaning that they satisfy the transversality properties
stated above and can therefore be used to define $HF_*(S^1,F^\pm)$.
These give rise to almost complex structures $J^\pm$ on the products
$\R \times F^\pm$, and the requirements are that $j$ is standard on
$(-\infty;s_-] \times S^1$ and $[s_+;\infty) \times S^1$, and
$\Psi^{\pm}$ is $(J,J^\pm)$-holomorphic. We denote the space of such
pairs $(j,J)$, for a fixed choice of $J^\pm$, by $\JJ(S,E)$.

For $\sigma_- \in \SS^h(S^1,F^-)$, $\sigma_+ \in \SS^h(S^1,F^+)$,
consider the space $\MM(S,E;\sigma_-,\sigma_+)$ of sections $u: S
\rightarrow E$ satisfying the same equation \eqref{eq:floer} as
before, where the convergence conditions should be more properly
formulated as $\Psi^\pm(u(s,t)) = (s,u^\pm(s,t))$ with
$u^\pm(s,\cdot) \rightarrow \sigma_\pm$ in $\SS(S^1,F^\pm)$.
Writing $\MM_0(S,E;\sigma_-,\sigma_+)$ for the subspace where the
virtual dimension is zero, one defines a chain homomorphism
$CG(S,E): CF_*(S^1,F^-) \rightarrow CF_*(S^1,F^+)$ by
\begin{equation} \label{eq:cg-map}
 CG(S,E)\gen{\sigma_-} = \!\!\!\! \sum_{
 \substack{
 \sigma_+ \in \SS^h(S^1,F^+) \\
 \!\!\!\! u \in \MM_0(S,E;\sigma_-,\sigma_+)\!\! } }
 \pm q^{\epsilon(u)} \gen{\sigma_+}
\end{equation}
The relative Gromov invariant is the induced map on homology. As
before, there are a number of conditions that $J$ has to satisfy, in
order for \eqref{eq:cg-map} to be a well-defined and meaningful
expression:
\begin{itemize}
\item
There are no $J$-holomorphic spheres in any fibre of $E$ with
strictly negative Chern number.
\item
If $v: S^2 \rightarrow E_{s,t}$ is a non-constant $J$-holomorphic
sphere with Chern number zero, its image does not contain $u(s,t)$
for any $u \in \MM_0(S,E;\sigma_-,\sigma_+)$.
\item
The linearized operator $D\bar\partial_u$ associated to any $u \in
\MM(S,E;\sigma_-,\sigma_+)$ is onto.
\end{itemize}
Note that because our fibration is a two-parameter family of
symplectic four-manifolds, pseudo-holomorphic spheres in the fibres
with Chern number zero can no longer be avoided, even though one can
always achieve that a particular fixed fibre contains none of them.
The proof that the above conditions are generic is standard; for
details consult \cite{hofer-salamon95} and \cite{mcduff-salamon}.

There is little difficulty in replacing our symplectic fibration
with a Lefschetz fibration $\pi: E \rightarrow S$, having the same
kind of behaviour at infinity. In this case, the definition of
$\JJ(S,E)$ includes the additional requirements that $j = j_S$ in
a neighbourhood of the critical values, and $J = J_E$ near the
critical points. A smooth section cannot pass through any critical
point, so the analytic setup for the moduli spaces
$\MM_0(S,E;\sigma_-,\sigma_+)$ remains the same as before. Of
course, pseudo-holomorphic spheres in the singular fibres appear
in the Gromov-Uhlenbeck compactification of the space of sections,
and to avoid potential problems with them one has to impose
another condition on $J$:
\begin{itemize}
\item
If $(s,t) \in S^{crit}$ and $v: S^2 \rightarrow E$ is a nonconstant
$J$-holomorphic map with image in $E_{s,t}$, then $\leftsc c_1(E),[v]
\rightsc > 0$.
\end{itemize}
To prove genericity of this, one considers the minimal resolution
$\hat{E}_{s,t}$ of $E_{s,t}$, which is well-defined because our
complex structure $J$ is integrable near the singularities. It is
a feature of ordinary double points in two complex dimensions
(closely related to simultaneous resolutions) that
$c_1(\hat{E}_{s,t})$ is the pullback of $c_1(E)|E_{s,t}$. By a
small perturbation of the almost complex structure on the
resolution, supported away from the exceptional divisor, one can
achieve that there are no pseudo-holomorphic curves $\hat{v}: S^2
\rightarrow \hat{E}_{s,t}$ with $\leftsc
c_1(\hat{E}_{s,t}),[\hat{v}] \rightsc \leq 0$ except for the
exceptional divisor itself and its multiple covers. The desired
result follows by lifting pseudoholomorphic spheres from $E_{s,t}$
to the resolution.

\subsection{}
Solutions of Floer's equation have two properties not shared by more
general pseudoholomorphic sections: (1) there is an $\R$-action by
translations; (2) the energy of any pseudoholomorphic section is
$\epsilon(u) \geq 0$, and those with zero energy are horizontal
sections of the symplectic connection. While (1) is characteristic of
Floer's equation, (2) can be extended to a wider class of geometric
situations, as follows. Let $\pi: E \rightarrow S = \R \times S^1$ be
a Lefschetz fibration with the same ``tubular end'' structure as
before. We say that $E$ has {\em nonnegative (Hamiltonian) curvature}
if for any point $x \notin E^{crit}$, the restriction of $\Omega$ to
$TE_x^h$ is nonnegative with respect to the orientation induced from
$TS_{\pi(x)}$. A pair $(j,J) \in \JJ(S,E)$ is {\em fully compatible}
if $\Omega(\cdot,J\cdot)$ is symmetric, or equivalently $J(TE^h_x)
\subset TE^h_x$ for all $x \notin E^{crit}$. With respect to the
decomposition \eqref{eq:jj} this means that $J^{vh} = 0$. The
following result is straightforward:

\begin{lemma} \label{th:nonnegative}
Suppose that $E$ has nonnegative curvature, and that $J$ is fully
compatible. Then any $u \in \MM(S,E;\sigma_-,\sigma_+)$ satisfies
$\epsilon(u) \geq 0$. Any horizontal (covariantly constant) section
is automatically $J$-holomorphic; in the converse direction, any $u
\in \MM(S,E;\sigma_-,\sigma_+)$ with $\epsilon(u) = 0$ must
necessarily be horizontal. \qed
\end{lemma}

To take advantage of this, one would like to make the spaces of
pseudo-holomor\-phic sections regular by choosing a generic $J$
within the class of fully compatible almost complex structures. It is
easy to see that all non-horizontal $u \in
\MM(S,E;\sigma_-,\sigma_+)$ can be made regular in this way, but the
horizontal sections persist for any choice of $J$, so we have to
enforce their regularity by making additional assumptions. The
following Lemma is useful for that purpose:

\begin{lemma}
In the situation of Lemma \ref{th:nonnegative}, let $u$ be a
horizontal section. Suppose that $\epsilon(u) = 0$, and that the
associated linearized operator $D\bar\partial_u$ has index zero. Then
$u$ is regular, which is to say that $D\bar\partial_u$ is onto. \qed
\end{lemma}

This is an easy consequence of a Weitzenb{\"o}ck argument, see
\cite[Lemmas 2.11 and 2.27]{seidel01}. Hence, if any horizontal $u$
satisfies the condition of the Lemma, one can indeed choose a fully
compatible $J$ which makes the moduli spaces of pseudo-holomorphic
sections regular. Full compatibility does not restrict the behaviour
of $J$ on the fibres, so we can also achieve all the other conditions
needed to make the relative Gromov invariant well-defined. After
expanding the resulting chain homomorphism into powers of $q$,
\begin{equation} \label{eq:cg}
 CG(S,E) = \sum_{d \geq 0} CG(S,E)_d q^d
\end{equation}
one finds that the leading term $CG(S,E)_0$ counts only horizontal
sections.

\begin{lemma} \label{th:nonneg}
Suppose that $E$ has nonnegative curvature, and that for any
$\sigma_+ \in \SS^h(S^1,F^+)$ there is a horizontal section $u$ of
$E$ with $\lim_{s \rightarrow +\infty} u(s,\cdot) = \sigma_+$, such
that $\epsilon(u) = 0$ and $D\bar\partial_u$ has index zero. Then,
for a suitable choice of fully compatible almost complex structure
$J$, the cochain level map $CG(S,E): CF_*(S^1,F^-) \rightarrow
CF_*(S^1,F^+)$ is surjective.
\end{lemma}

\proof Since horizontal sections are determined by their value at any
single point, it follows that for any $\sigma_+$ there is a unique
horizontal section $u = u_{\sigma_+}$ approaching it, and that these
are all horizontal sections. Consider the map $R: CF_*(S^1,F^+)
\rightarrow CF_*(S^1,F^-)$ which maps $\gen{\sigma_+}$ to the
generator $\gen{\sigma_-}$ associated to the negative limit of
$u_{\sigma_+}$. From \eqref{eq:cg} one sees that $(CG(S,E) \circ
R)\gen{\sigma_+} = \pm\gen{\sigma_+} + \text{(strictly positive
powers of $q$)}$, which clearly shows that $CG(S,E)$ is onto. \qed

We now turn to the concrete problem posed by a Lagrangian sphere $L$
in a symplectic four-manifold $M$. Choose a symplectic embedding $i:
T^*_{\scriptscriptstyle <\lambda}S^2 \rightarrow M$ with $i(S^2) =
L$. Take a model Dehn twist $\tau$ defined using a function $r$ which
satisfies
\begin{equation} \label{eq:wobbly}
 r'(t) \begin{cases} \in [1/4;3/4] & t \in [0;\mu), \\
                     \in [0;1/2] & t \in [\mu;\lambda/2), \\
                     = 0 & t \geq \lambda/2
\end{cases}
\end{equation}
for some $\mu < \lambda/2$, and transplant it to a Dehn twist
$\tau_L$ using $i$. Next, choose a a Morse function $H: M \rightarrow
\R$ with the properties that (1) $H(i(u,v)) = ||u||$ for all $\mu
\leq ||u|| < \lambda$; (2) $h$ has precisely two critical points in
$im(i)$, both of which lie on $L$ (their Morse indices will of course
be 0 and 2). Let $\phi$ be the Hamiltonian flow of $H$ for small
positive time $\delta>0$, and $\tilde{\tau}_L = \tau_L \circ \phi$.
By construction, $\tilde{\tau}_L = \phi$ outside
$i(T^*_{\scriptscriptstyle <\lambda/2}S^2)$.

\begin{lemma} \label{th:no-fixed-points}
If $\delta$ is sufficiently small, $\tilde{\tau}_L$ has no fixed
points inside $im(i)$.
\end{lemma}

\proof As $\delta$ becomes small, the fixed points of
$\tilde{\tau}_L$ accumulate at the fixed points of $\tau_L$, hence
they will lie outside $i(T^*_{\scriptscriptstyle \leq \mu} S^2)$.
Recalling the definition of the model Dehn twist, we have that
\[
 (i^{-1}\tilde{\tau}_L i)(u,v) = \sigma_{2\pi (r'(||u||) + \delta)}(u,v)
\]
for $||u|| \geq \mu$, and since $r'(||u||) + \delta \in
[\delta;1/2+\delta]$ cannot be an integer, $i(u,v)$ cannot be a
fixed point. \qed

Let $S' \subset S = \R \times S^1$ be the disc of radius $1/4$ around
$(s,t) = (0,0)$, and $S'' = S \setminus int(S')$. There is a
Lefschetz fibration $E' \rightarrow S'$ with fibre $M$, whose
monodromy around $\partial D$ is a Dehn twist $\tau_L$ defined using
a function $r$ that satisfies \eqref{eq:wobbly}. This is explicitly
constructed in \cite[Section 1.2]{seidel01}, where properties
somewhat stricter to \eqref{eq:wobbly} are subsumed under the notion
of ``wobblyness''. To complement this, there is a fibration $E''
\rightarrow S''$ with a two-form $\Omega''$ such that the monodromy
of the resulting symplectic connection is $\phi$ around the loop
$\{-1\} \times S^1$, $\tau_L \circ \phi = \tilde{\tau}_L$ around
$\{+1\} \times S^1$, and $\tau_L$ around $\partial S''$. This is
actually much simpler to write down:
\[
 E'' = \frac{\{z = s+it \in \R \times [0;1] \suchthat
 |z| \geq 1/4, |z-i| \geq 1/4\}}{(s,1,x) \sim (s,0,\phi(x))
 \text{ for $s < 0$}, \; (s,1,x) \sim (s,0,\tilde{\tau}_L(x))
 \text{ for $s > 0$}}
\]
One can glue together the two pieces along $\partial S' = \partial
S''$ to a Lefschetz fibration $E \rightarrow S$, and the resulting
chain level map is
\begin{equation} \label{eq:cg-final}
 CG(S,E): CF_*(\phi) \longrightarrow CF_*(\tilde{\tau}_L).
\end{equation}
It can be arranged that $E'$ has nonnegative curvature, and that
$E''$ is flat (zero curvature), so the curvature of $E$ is again
nonnegative. Actually, the construction of $E'$, like that of
$\tau_L$ itself, is based on the local model of
$T^*_{\scriptscriptstyle \leq \lambda} S^2$, so that $E'$ contains a
trivial piece $S' \times (M \setminus im(i))$. Using this and Lemma
\ref{th:no-fixed-points} one sees that for any fixed point $x$ of
$\tilde{\tau}_L$, which the same as a critical point of $H$ lying
outside $im(i)$, there is a horizontal section $u$ of $E$ such that
$u(s,t) = (s,t,x)$ for $s > 1/2$, and that these sections satisfy the
conditions of Lemma \ref{th:nonneg}. Hence \eqref{eq:cg-final} is
onto for a suitable choice of almost complex structure; but from the
definition of $H$, we know that its kernel is two-dimensional and
concentrated in $CF_{even}(\phi)$, which implies that the induced map
$G = G(S,E)$ on Floer homology fits into a long exact sequence as
stated in Proposition \ref{th:exact}. On the other hand, a gluing
argument which separates the two pieces in our construction of $E$
shows that one can indeed write $G$ as pair-of-pants product with an
element $\theta_L$ as in \eqref{eq:char}.

\providecommand{\bysame}{\leavevmode\hbox
to3em{\hrulefill}\thinspace}
\providecommand{\MR}{\relax\ifhmode\unskip\space\fi MR }
\providecommand{\MRhref}[2]{%
  \href{http://www.ams.org/mathscinet-getitem?mr=#1}{#2}
} \providecommand{\href}[2]{#2}

\end{document}